# HEAVY-TRAFFIC LIMITS FOR WAITING TIMES IN MANY-SERVER QUEUES WITH ABANDONMENT[1]

By Rishi Talreja and Ward Whitt

*Columbia University*

We establish heavy-traffic stochastic-process limits for waiting times in many-server queues with customer abandonment. If the system is asymptotically critically loaded, as in the quality-and-efficiency-driven (QED) regime, then a bounding argument shows that the abandonment does not affect waiting-time processes. If instead the system is overloaded, as in the efficiency-driven (ED) regime, following Mandelbaum et al. [*Proceedings of the Thirty-Seventh Annual Allerton Conference on Communication, Control and Computing* (1999) 1095–1104], we treat customer abandonment by studying the limiting behavior of the queueing models with arrivals turned off at some time $t$. Then, the waiting time of an infinitely patient customer arriving at time $t$ is the additional time it takes for the queue to empty. To prove stochastic-process limits for virtual waiting times, we establish a two-parameter version of Puhalskii's invariance principle for first passage times. That, in turn, involves proving that two-parameter versions of the composition and inverse mappings appropriately preserve convergence.

**1. Introduction.** In the past, heavy-traffic stochastic-process limits for waiting times in queueing systems have been proven using Puhalskii's invariance principle for first-passage times [11] together with established stochastic-process limits for queue-length processes; e.g., [1, 4, 13, 19]. For instance, for an $n$ server system, denoting the arrival, departure and number-in-system processes by $A$, $D$ and $X$, respectively, the virtual waiting time at time $t$ can be represented as the following first-passage time:

$$(1.1) \qquad V(t) = \inf\{s \geq 0 \mid D(t+s) \geq X(0) + A(t) - (n-1)\}.$$

Received April 2008; revised March 2009.

[1]Supported by NSF Grant DMI-0457095.

When one attempts to incorporate an abandonment process $L$ in the model, one cannot simply write the virtual waiting time at time $t$ (the waiting time of an infinitely patient customer arriving at time $t$) as

$$(1.2) \quad V(t) = \inf\{s \geq 0 \mid D(t+s) + L(t+s) \geq X(0) + A(t) - (n-1)\},$$

because the term $L(t+s)$ may include customers that arrived to the system after time $t$ and then abandoned.

When the system is asymptotically critically loaded in fluid scale, as in the quality-and-efficiency (QED) regime, we show that abandonments do not affect fluid and diffusion limits for the virtual waiting time. We do this by bounding the virtual waiting time process from above and below by processes of the form (1.1) and (1.2), respectively. These bounds are shown to be asymptotically equivalent, giving us a limit for our virtual waiting time process. For the Markovian $M/M/n + M$ model in the QED regime, the virtual waiting time limit was first stated in Theorem 3 of Garnett et al. [3], but the proof given there is not correct. We give a correct proof here.

When the system is not asymptotically critically loaded in fluid scale, we must be more careful. One way to address this problem, which was first suggested by Mandelbaum et al. [7], is to consider the system with the arrival process stopped at a fixed time $t$, so that the abandonment term in the first-passage time expression (1.2) does not include customers that arrive to the system after time $t$. We can then use the invariance principle for first-passage times to get the desired limit for the fixed $t$. To prove the full stochastic-process limit for virtual waiting times, however, we need a corresponding two-parameter limit for the queue-length process with arrivals turned off and a two-parameter version of the invariance principle for first-passage times. We establish these results here. We establish the two-parameter version of Puhalskii's invariance principle for first-passage times by proving that two-parameter versions of the composition and inverse mappings are convergence-preserving. This portion of our work is along the lines of [14] and Chapter 13 of [15].

We apply the two-parameter version of Puhalskii's theorem to establish many-server heavy-traffic stochastic-process limits for waiting times in the $M/M/s + M$ model in the efficiency-driven (ED) regime. As indicated above, the line of reasoning here was initiated by Mandelbaum et al. [7], which builds on the strong-approximation approach in [6]. In [7] the need to turn off the arrival process to properly treat customer abandonment was recognized, but the two-parameter version of Puhalskii's theorem needed to complete the argument was not stated or proved. The statement and partial proof were developed in an unpublished manuscript [8]. The most difficult step in proving the ED result is not establishing the two-parameter generalization of Puhalskii's theorem, but is instead establishing the stochastic-process limit for the queue-length process with the arrival process turned off, Proposition 6.2 here. Mandelbaum et al. [8] establish corresponding results in the more general setting of "Markovian service networks" [6], but that part of their work remains to be completed. Our task is made easier by considering only the $M/M/s + M$ model. Nevertheless, we believe that the results and line of reasoning here apply to more general models.

We also establish the heavy-traffic limit for the steady-state waiting time in the ED regime by a different argument. As should be anticipated, it is a normal distribution centered at the fluid limit. We show that the variance has a simple form, which is unclear a priori. Our result is a special case of the ED limit for the $M/M/s + GI$ steady-state waiting time in (6.12) of Zeltyn and Mandelbaum [20]. Our proof is interesting because we use a stochastic-process limit to prove the one-dimensional limit. For the QED regime, the limit for the steady-state waiting time distribution is given in Theorem 3 of [3].

NOTATION.  For $a, b, c \in \mathbb{R}$, let $a \vee b \equiv \max\{a, b\}$, $a \wedge b \equiv \min\{a, b\}$, $[c]_a^b \equiv \max\{\min\{c, b\}, a\}$, and $a^+ \equiv a \vee 0$.

For a complete separable metric space (CSMS) $(\mathcal{S}, d)$, we let $C([0, \infty), \mathcal{S})$ denote the set of continuous functions from $[0, \infty)$ to $\mathcal{S}$ and define $C \equiv C([0, \infty), \mathbb{R})$ and $C_C \equiv C([0, \infty), C)$. We endow these spaces with the topology of uniform convergence over compact sets (u.o.c.). For $x \in C([0, \infty), \mathcal{S})$, let $\|x - y\|_t \equiv \sup_{0 \leq s \leq t} d(x(s), y(s))$. A metric $m$ inducing the topology of u.o.c. is

$$(1.3) \qquad m(x, y) \equiv \int_0^\infty e^{-t}(\|x - y\|_t \wedge 1)\, dt.$$

Notice that $m(x_n, x) \to 0$ as $n \to \infty$ if and only if $\|x_n - x\|_t \to 0$ as $n \to \infty$ for all $t \geq 0$.

For a CSMS $\mathcal{S}$, let $D([0, \infty), \mathcal{S})$ denote the set of right-continuous functions from $[0, \infty)$ to $\mathcal{S}$ having left limits everywhere on $(0, \infty)$ (see [2] and [15]) and let $D \equiv D([0, \infty), \mathbb{R})$ and $D_D \equiv D([0, \infty), D)$. Let $e$ denote the identity element in $D$, i.e., $e(t) = t$ for $t \geq 0$, and let $\eta$ denote the constant function 1, i.e., $\eta(t) = 1$ for $t \geq 0$. We use $(D, \mathcal{T})$ to denote the space $D$ with the topology $\mathcal{T}$. We use $(D_D, \mathcal{T}_1, \mathcal{T}_2)$ to denote the space $D_D$ with the topology $\mathcal{T}_1$ on the outside $D$ space and the topology $\mathcal{T}_2$ on the inside $D$ space. We use $J_1$ and $M_1$ to denote Skohorod's $J_1$ and $M_1$ topologies, respectively, and $U$ to denote the topology of u.o.c. as defined above. Note the spaces $(D_D, J_1, M_1)$ and $(D_D, J_1, J_1)$ are well-defined topological spaces, but $(D_D, M_1, M_1)$ and $(D_D, M_1, J_1)$ are not, since the $M_1$ topology has only been defined on the spaces $D([0, \infty), \mathbb{R}^k)$, for $k \geq 1$ (see Chapter 12 of [15]).

We can identify functions in $C_C$ and $D_D$ as two-parameter functions $x \equiv x(\cdot, \cdot)$. For $x \in D_D$ and $s \geq 0$, we will often find it convenient to refer to $x(s, \cdot) \in D([0, \infty), \mathbb{R})$ as $x^s$. We also work with product spaces such as $D_D^2 \equiv D_D \times D_D$, where we assume the appropriate product topology. For $x \in D$, let $\|x\|_S \equiv \sup_{0 \leq s \leq S} |x(s)|$ and for $x \in D_D$, let $\|x\|_{S,T} \equiv \sup_{0 \leq s \leq S, 0 \leq t \leq T} |x(s, t)|$. (This is consistent with our notation $\|x - y\|_t$ used in (1.3); since $\mathcal{S} = \mathbb{R}$, we now have a norm.) We will mainly be considering limit processes with continuous sample paths. It is thus significant that for $x \in C_C$ convergence $x_n \to x$ in $(D_D, J_1, J_1)$ or $(D_D, J_1, M_1)$ is equivalent to $x_n \to x$ in $(D_D, U, U)$, which, in turn, is equivalent to $\|x_n - x\|_{S,T} \to 0$ for all $S, T > 0$; see Theorem 2.1 below.

Let $\Rightarrow$ denote convergence in distribution. For a sequence of $\mathbb{R}$-valued random variables $(X_n)_{n \geq 1}$, we write $X_n = o_{\mathbb{P}}(f(n))$ if

$$\frac{X_n}{f(n)} \Rightarrow 0 \qquad \text{in } \mathbb{R} \text{ as } n \to \infty.$$

*Organization.* In Section 2, we begin by proving various functions in $D_D$ are convergence preserving. Some of these results are given in greater generality than needed here. The additional generality may be useful in future work. In Section 3, we use a simple bounding argument to prove heavy-traffic limits for virtual waiting time processes when the abandonment process is asymptotically negligible in fluid-scale. In Section 4, we use the results in Section 2 to prove heavy-traffic limits for virtual waiting time processes assuming associated limits hold for basic queueing process with the arrival process stopped. In Sections 5 and 6 we use the results in Sections 3 and 4, respectively, to prove waiting-time stochastic-process limits for the $M/M/n + M$ model in the QED and ED limiting regimes. In Section 6.2, we establish the heavy-traffic limit for the steady-state virtual waiting times in the ED regime. Section 7 contains the proofs of lemmas used in Sections 3 and 6.

**2. Functions in $D_D$.** In this section, we prove continuity results for functions in $D_D$. We focus on composition and inverse in Sections 2.1 and 2.2, respectively. We use the results in Sections 2.1 and 2.2 to extend Puhalskii's invariance principle for first-passage times [11] to functions of two parameters in Section 2.3. We prove continuity results for integral mappings in Section 2.4 and for a projection mapping in Section 2.5. For our two-parameter results, the space can be either $(D_D, J_1, J_1)$ or $(D_D, J_1, M_1)$, unless the $M_1$ topology is explicitly assumed on the inside space, but this usually does not matter because our limits will usually be assumed to be in $C_C$.

Let $D_{2,0}$ be the subset of functions $x$ in $D_D$ such that $x(0, 0) \geq 0$. Let $D_{2,\uparrow}$ and $D_{2,\uparrow\uparrow}$ be the subsets of functions in $D_{2,0}$ that are nondecreasing and strictly increasing in each coordinate, respectively. Similarly, define $C_{2,0} \equiv C_C \cap D_{2,0}$, $C_{2,\uparrow} \equiv C_C \cap D_{2,\uparrow}$ and $C_{2,\uparrow\uparrow} \equiv C_C \cap D_{2,\uparrow\uparrow}$.

2.1. *Composition.* First, we show that, as in the one-parameter case, convergence in $D_D$ to a point in $C_C$ is equivalent to uniform convergence over compact sets. We need the following lemma, whose proof is an easy exercise in analysis.

LEMMA 2.1. *Let $(\mathcal{S}, d)$ be a metric space. Let $(x_n)_{n \geq 1}$ and $x$ be functions from a compact set $E \subset \mathbb{R}^k$ to $\mathcal{S}$, $k \geq 1$, and let $x$ be continuous. Then $\sup_{t \in E} d(x_n(t), x(t)) \to 0$ as $n \to \infty$ if and only if for all sequences $(t_n)_{n \geq 1} \subseteq E$ such that $t_n \to t$ as $n \to \infty$ we have $x_n(t_n) \to x(t)$ in $\mathcal{S}$ as $n \to \infty$.*

THEOREM 2.1. *If $x \in C_C$, then $x_n \to x$ in $(D_D, J_1, J_1)$ or $(D_D, J_1, M_1)$ as $n \to \infty$ if and only if $x_n \to x$ in $(D_D, U, U)$ as $n \to \infty$.*

PROOF. Let $\mathcal{T}$ be $J_1$ or $M_1$. Since $x \in C_C$, the convergence $x_n \to x$ in $(D_D, J_1, \mathcal{T})$ is equivalent to convergence in $(D_D, U, \mathcal{T})$. Therefore, for all $T > 0$, we have

$$\tag{2.1} \sup_{t \in [0,T]} d_{\mathcal{T}}(x_n^t, x^t) \to 0 \qquad \text{as } n \to \infty,$$

where $d_{\mathcal{T}}$ is a metric on $D$ inducing the topology $\mathcal{T}$. By Lemma 2.1, (2.1) holds if and only if for every $T > 0$, $t \in [0, T]$, and sequence $(t_n)_{n \geq 1} \subseteq [0, T]$ such that $t_n \to t$ as $n \to \infty$ we have $x_n^{t_n} \to x^t$ in $(D([0, T], \mathbb{R}), \mathcal{T})$. But since $x^t \in C$ for each $t \in [0, T]$ this is equivalent to $x_n^{t_n} \to x^t$ in $(D([0, T], \mathbb{R}), U)$ for all $T > 0$. Again, applying Lemma 2.1, this convergence holds if and only if for for each $S, T > 0$, $(t, s) \in [0, T] \times [0, S]$, and sequence $((t_n, s_n))_{n \geq 1} \subseteq [0, T] \times [0, S]$ such that $(t_n, s_n) \to (t, s)$ as $n \to \infty$ we have $x_n(t_n, s_n) \to x(t, s)$ in $\mathbb{R}$ as $n \to \infty$. But then again by Lemma 2.1, this is true if and only if $x_n \to x$ in $(D_D, U, U)$.  □

We now show that addition and multiplication on $D_D^2$ are measurable and continuous at points in $C_C^2$. In both cases, to prove measurability, we use Lemma 2.7 of [14], which we state here for completeness.

LEMMA 2.2 (Lemma 2.7 of [14]). *The Borel $\sigma$-field on $D_D$ coincides with the Kolmogorov $\sigma$-field, i.e., the $\sigma$-field generated by the coordinate projections.*

Let addition on $D_D^2$ be defined pointwise: For $x, y \in D_D$,

$$(x + y)(s, t) \equiv x(s, t) + y(s, t).$$

THEOREM 2.2 (Continuity of addition in $D_D$). *Addition on $D_D^2$ is measurable and continuous at points in $C_C^2$.*

PROOF. Measurability follows by Theorem 4.1 of [14]. The argument there is based on Lemma 2.2. Continuity follows from Theorem 2.1 and the triangle inequality. For each $S, T > 0$,

$$\|(x_n + y_n) - (x + y)\|_{S,T} \leq \|x_n - x\|_{S,T} + \|y_n - y\|_{S,T}. \qquad \square$$

Similarly, define multiplication on $D_D$ by

$$(xy)(s,t) \equiv x(s,t)y(s,t)$$

for $x, y \in D_D$.

THEOREM 2.3 (Continuity of multiplication in $D_D$). *Multiplication on $D_D^2$ is measurable and continuous at points in $C_C^2$.*

PROOF. Again, measurability follows using the same argument as in Theorem 4.1 of [14]. Continuity at continuous limits follows from Theorem 2.1 and the triangle inequality:

$$\|x_n y_n - xy\|_{S,T} \leq \|x_n y_n - x_n y\|_{S,T} + \|x_n y - xy\|_{S,T}$$

$$\leq \|x_n\|_{S,T}\|y_n - y\|_{S,T} + \|y\|_{S,T}\|x_n - x\|_{S,T}.$$

Since $x_n \to x$, $\sup_{n \geq 1} \|x_n\|_{S,T} < \infty$, and we have our result. $\qquad \square$

Define the composition mapping $\circ_2$ from $D_D \times D_D$ to $D_{2,\uparrow}$ by

$$(2.2) \qquad (x \circ_2 y)(s,t) \equiv x(s, y(s,t)).$$

Notice that $y$ only enters in the second argument of $x$. As usual, we require that the internal function be in $D_{2,\uparrow}$. This is needed to ensure that the composition $x \circ_2 y$ is an element of $D_D$ (see Example 13.2.1 of [15]). Alternatively, we can define the composition mapping as a mapping from $C_C \times C_C$ to $C_C$. For each $s \geq 0$, we have $(x \circ_2 y)(s, \cdot) = x(s, y(s, \cdot)) = x^s \circ y^s$, where $\circ$ is the standard composition map on $D$.

Then we have the following extension of Theorem 3.1 of [14] and Theorem 13.2.2 of [15].

THEOREM 2.4 (Convergence preservation of $\circ_2$).

1. *The composition function $\circ_2$ mapping $D_D \times D_{2,\uparrow}$ into $D_D$ defined in (2.2) is measurable and continuous at each $(x,y) \in C_C \times C_{2,\uparrow}$.*
2. *The composition function $\circ_2$ mapping $C_C \times C_C$ into $C_C$ is continuous.*

PROOF. For part (1), measurability follows from Lemma 2.2 and the fact that the composition mapping from $D \times D_\uparrow$ to $D$ is measurable (see

Theorem 13.2.1 of [15] and page 232 of [2]). For both parts (1) and (2), we use the triangle inequality to prove continuity. For each $S, T > 0$, we have

$$\|x_n \circ_2 y_n - x \circ_2 y\|_{S,T} \leq \|x_n \circ_2 y_n - x \circ_2 y_n\|_{S,T} + \|x \circ_2 y_n - x \circ_2 y\|_{S,T}.$$

Let $T' \equiv \sup_n \|y_n\|_{S,T}$. Then the first term converges to zero since it is bounded by $\|x_n - x\|_{S,T'} \to 0$. The second term converges to zero since $x \in C_C$. The mapping in part (2) is measurable because it is continuous. $\square$

We now give the main result of this subsection. It is an extension of Theorem 13.3.3 of [15].

THEOREM 2.5 (Convergence preservation of $\circ_2$ with nonlinear centering). *Let* $(x_n, y_n) \in (D_D \times D_{2,\uparrow}) \cup (C_C \times C_C)$ *and* $y \in C_C$. *Let* $x$ *have continuous partial derivative* $x'$ *with respect to its second parameter, i.e.* $x'(s, t) \equiv \frac{\partial x(s,t)}{\partial t}$, *and let* $c_n \to \infty$. *If*

$$c_n(x_n - x, y_n - y) \to (u, v) \qquad in \ D_D^2,$$

*where* $u, v \in C_C$, *then*

$$c_n(x_n \circ_2 y_n - x \circ_2 y) \to u \circ_2 y + (x' \circ_2 y)v \qquad in \ D_D.$$

PROOF. We follow the proof of Theorem 13.3.3 of [15]. Let $(x_n, y_n) \in D_D \times D_{2,\uparrow}$; the other case can be proven analogously. Note that

$$c_n(x_n \circ_2 y_n - x \circ_2 y) = c_n(x_n - x) \circ_2 y_n + c_n(x \circ_2 y_n - x \circ_2 y).$$

By our assumptions, we have

$$(c_n(x_n - x), c_n(y_n - y), y_n) \to (u, v, y) \qquad in \ D_D^3,$$

so, applying Theorems 2.2 and 2.4, gives us

$$(2.3) \quad (c_n(x_n \circ_2 y_n - x \circ_2 y_n) + (x' \circ_2 y)c_n(y_n - y)) \to u \circ_2 y + (x' \circ_2 y)v.$$

Now, for all $S, T > 0$,

$$\|c_n(x_n \circ_2 y_n - x \circ_2 y) - c_n(x_n \circ_2 y_n - x \circ_2 y_n) - c_n(x' \circ_2 y)(y_n - y)\|_{S,T}$$

$$(2.4) \quad \begin{aligned} &\leq \|c_n(x \circ_2 y_n - x \circ_2 y) - c_n(x' \circ_2 y)(y_n - y)\|_{S,T} \\ &= \sup_{\substack{0 \leq s \leq S \\ 0 \leq t \leq T}} \left| c_n \int_{y(s,t)}^{y_n(s,t)} x'(s, u) \, du - c_n x'(s, y(s,t))(y_n(s,t) - y(s,t)) \right| \\ &\leq \sup_{\substack{0 \leq s \leq S \\ 0 \leq t \leq T}} \left| \left( \sup_{y(s,t) \leq u \leq y_n(s,t)} x'(s, u) - x'(s, y(s,t)) \right) c_n(y_n(s,t) - y(s,t)) \right|. \end{aligned}$$

The last expression goes to zero since $x'$ is uniformly continuous over bounded intervals in its first parameter and $c_n(y_n - y) \to v$ in $D_D$. Now, combining (2.3) and (2.4) and using the triangle inequality gives us our result. $\square$

2.2. *Inverse with centering.* We now move on to prove results about inverse maps with linear and nonlinear centering in $D_D$. The inverse we consider here is the usual inverse applied to only the second argument. Let $D_{2,u}$ be the subset of functions $x$ in $D_D$ such that $x(0,0) = 0$ and $x(s,\cdot)$ is unbounded above for each $s \geq 0$. As before, let $D_{2,u,\uparrow} \equiv D_{2,u} \cap D_{2,\uparrow}$ and $D_{2,u,\uparrow\uparrow} \equiv D_{2,u} \cap D_{2,\uparrow\uparrow}$. We define the inverse map on the subset $D_{2,u}$ of $D$ as follows. For $x \in D_{2,u}$, let the inverse of $x$ be

$$x^{-1}(s,t) \equiv \inf\{u \geq 0 \mid x(s,u) > t\}.$$

Notice that for each $s \geq 0$,

$$(x^{-1})^s = x^{-1}(s,\cdot) = \inf\{u \mid x(s,u) > \cdot\} = \inf\{u \mid x^s(u) > \cdot\} = (x^s)^{-1},$$

where the inverse on the right is the standard inverse on $D$.

We then have the following result, which is analogous to Theorem 13.6.1 of [15], which is proved in Section 7.6.1 of [16]. The theorem exploits the $M_1$ topology on the inner space $D$. Note the second condition of the theorem holds if $x \in C_C$, which will be the case when we apply the theorem below.

THEOREM 2.6 (Convergence preservation of inverse). *The inverse map on $(D_{2,u}, J_1, M_1)$ is measurable and continuous at all $x \in D_{2,u}$ such that:*

1. $x^{-1}(s,0) = 0$ *for all $s \geq 0$,*
2. $x$ *is right-continuous at 0, uniformly in $s \in [0,S]$ for each $S > 0$, i.e., for all $S > 0$ and $\varepsilon > 0$, there exists $\delta > 0$ such that*

$$(2.5) \qquad\qquad \sup_{0 \leq s \leq S} |x^s(\delta) - x^s(0)| < \varepsilon.$$

PROOF. Measurability follows by Lemma 2.2 and the fact that the inverse map on $D$ is measurable, by the argument in the proof of Theorem 13.6.1 of [15].

Suppose $x_n \to x$ in $(D_{2,u}, J_1, M_1)$. Then, for each $S > 0$, there exists a sequence of increasing homeomorphisms on $[0,S]$ such that

$$\|\lambda_n - e\|_S \vee \sup_{0 \leq s \leq S} d_{M_1}(x_n^s, x^{\lambda_n(s)}) \to 0 \qquad \text{as } n \to \infty.$$

It now suffices to show that $(x_n^{-1})^s \to (x^{-1})^{\lambda_n(s)}$ in $(D, M_1)$, uniformly in $s \in [0,S]$. To show this, we can use the same continuity argument as in the proof of Theorem 13.6.1 of [15]. In particular, we can construct a sequence of functions $x_n^*$ in $D_D$ such that for $n$ large enough $(x_n^*)^{-1}(s,0) = 0$ for all $s \in [0,S]$ and $((x_n^*)^s)^{-1}$ is close to $(x_n^s)^{-1}$ in $(D,U)$ uniformly in $s \in [0,S]$. In order for that argument to go through uniformly in $s \in [0,S]$, we need the condition (2.5). The reason we want $(x_n^*)^{-1}(s,0) = 0$ for each $s \in [0,S]$ is that if we have this property then each $M_1$ parametric representation of $(x_n^*)^s$ serves as a parametric representation of $((x_n^*)^s)^{-1}$ with the roles of the coordinates switched (see Lemma 13.6.8 of [15]). Therefore, for large enough $n$,

$$\sup_{0 \le s \le S} d_{M_1}((x_n^{-1})^s, (x^{-1})^{\lambda_n(s)})$$

$$= \sup_{0 \le s \le S} d_{M_1}((x_n^s)^{-1}, (x^{\lambda_n(s)})^{-1})$$

$$\le \sup_{0 \le s \le S} d_{M_1}(((x_n^*)^s)^{-1}, (x^{\lambda_n(s)})^{-1}) + \sup_{0 \le s \le S} d_{M_1}(((x_n^*)^{-1})^s, (x_n^s)^{-1})$$

$$\le \sup_{0 \le s \le S} d_{M_1}(((x_n^*)^s)^{-1}, (x^{\lambda_n(s)})^{-1}) + \sup_{0 \le s \le S} \|((x_n^*)^{-1})^s - (x_n^s)^{-1}\|$$

$$= \sup_{0 \le s \le S} d_{M_1}((x_n^*)^s, x^{\lambda_n(s)}) + \sup_{0 \le s \le S} \|((x_n^*)^{-1})^s - (x_n^s)^{-1}\| \to 0$$

as $n \to \infty$. $\square$

REMARK 2.1. Theorem 7.1 of [14], which proves continuity of the inverse function in $(D, M_1)$, is incorrect. We need the limit process $x$ to have the extra property that $x^{-1}(0) = 0$. This has been corrected in Section 13.6 of [15].

We now prove convergence results for inverse with centering as in Section 13.7 of [15]. Here, we let the function $e_2 \in D_D$ be defined as $e_2(s, t) = t$ for all $s, t \ge 0$. Also, for each $T \ge 0$ define the maximum jump function $J_T : D \to \mathbb{R}$ by

$$J_T(x) \equiv \sup_{0 \le t \le T} \{|x(t) - x(t-)|\}.$$

THEOREM 2.7 (Inverse with linear centering). *Suppose that $c_n(x_n - e_2) \to x$ in $D_D$ as $n \to \infty$, where $x_n \in D_{2,u}$, $x \in C_C$, $y(0) = 0$, and $c_n \to \infty$. Then $c_n(x_n^{-1} - e_2) \to -x$ in $D_D$ as $n \to \infty$.*

PROOF. By the triangle inequality, for each $S, T > 0$,

$$\|c_n(x_n^{-1} - e_2) + x\|_{S,T} \le \|c_n(x_n^{-1} - x_n \circ_2 x_n^{-1}) + x\|_{S,T} + \|c_n(x_n \circ_2 x_n^{-1} - e_2)\|_{S,T}.$$

Applying the composition map, $\circ_2$ (Theorem 2.4) with the assumed convergence and $x_n^{-1} \to e_2$ (Theorem 2.6) gives us convergence of the first term above to 0. Furthermore, by Corollary 13.6.2 of [15], we have

$$\|c_n(x_n \circ_2 x_n^{-1} - e_2)\|_{S,T} = c_n \sup_{0 \le s \le S} \|x_n^s \circ (x_n^s)^{-1} - e\|_T$$

$$\le c_n \sup_{0 \le s \le S} J_{(x_n^s)^{-1}(T)}(x_n^s)$$

$$= \sup_{0 \le s \le S} J_{(x_n^s)^{-1}(T)}(c_n(x_n^s - e)).$$

Since $c_n(x_n - e) \to x$ and $x$ is in $C_C$, the last expression above converges to 0 as $n \to \infty$. This gives us our result. □

We now move on to inverse with nonlinear centering. The proofs of our next two results, Theorems 2.8 and 2.9, follow the proof of Theorems 13.7.2 and 13.7.4 of [15] exactly, with results on composition and inverse maps in $D$ there replaced by our new results on composition and inverse maps in $D_D$ here.

THEOREM 2.8 (Inverse with nonlinear centering). *Suppose that* $c_n(x_n - \lambda) \to u$ *as* $n \to \infty$ *in* $D_D$, *where* $x_n \in D_{2,u}$, $u \in C_C$, $u(0) = 0$, $c_n \to \infty$, *and* $\lambda$ *is absolutely continuous with continuous positive derivative* $\lambda'$ *with respect to its second parameter. Then*

$$c_n(x_n^{-1} - \lambda^{-1}) \to \frac{-u \circ_2 \lambda^{-1}}{\lambda' \circ_2 \lambda^{-1}} \qquad \text{in } D_D \text{ as } n \to \infty,$$

*where* $(u/v)(s,t) = u(s,t)/v(s,t)$, *for* $s, t \geq 0$.

2.3. *Two-parameter version of Puhalskii's theorem.* Finally, we have our two-parameter version of Puhalskii's theorem. This is similar to Theorem 1 in [11] and Theorem 13.7.4 in [15], which is an extension to the $M_1$ topology. Note that the limits $x$, $y$, $u$ and $v$ below are all continuous. Again, we omit the proof because it is a direct analog of the proof of Theorem 13.7.1 in [15], using the results above.

THEOREM 2.9. *Suppose that* $x_n \in D_{2,u}$, $y_n \in D_{2,\uparrow}$, $c_n \to \infty$,

$$c_n(x_n - x, y_n - y) \to (u, v) \qquad \text{in } D_D^2,$$

*where* $u, v \in C_C$, $x \in C_{2,u}$ *and* $y \in C_{2,\uparrow}$. *Furthermore, suppose* $x$ *is absolutely continuous with a continuous positive partial derivative* $x'$ *with respect to its second parameter, where* $x' \in C_C$. *Then*

$$(2.6) \quad c_n(x_n^{-1} \circ_2 y_n - x^{-1} \circ_2 y) \to \frac{v - u \circ_2 x^{-1} \circ_2 y}{x' \circ_2 x^{-1} \circ_2 y} \qquad \text{in } D_D \text{ as } n \to \infty.$$

REMARK 2.2. Although we require here that $x_n^s$ is unbounded for each $s \geq 0$ and $n \geq 1$, the proof goes through as long as the expressions on the left-hand side of (2.6) exist as elements of $D_D$. A sufficient conditions for this is $\sup_{u \geq 0} x_n^s(u) > \sup_{u \geq 0} y_n^s(u)$ for each $s \geq 0$ and $n \geq 1$.

REMARK 2.3. If we had defined the inverse process for $x \in D_{2,u}$ and $y \in D_{2,\uparrow}$ by

$$(x^{-1} \circ_2 y)(s,t) \equiv \inf\{u \geq 0 \mid x(s,u) \geq y(s,t)\},$$

then $x^{-1} \circ_2 y$ would not necessarily be an element of $D_D$. However, if $y^s$ is piecewise constant for each $s \geq 0$, then $x^{-1} \circ_2 y$ is an element of $D_D$. In fact, if we make the assumption that $y_n^s$ is piecewise constant for all $s \geq 0$ and $n \geq 1$, then the theorem holds with this alternative definition of inverse (see the remark in [11]). This will be useful in Section 4.

2.4. *Integral mappings.* We now prove continuity and measurability results for integral mappings in $D_D$, which are used to prove the two-parameter version of our ED queue-length heavy-traffic limit below. We first state a continuity result for the basic integral mapping for functions in $D_D$, where the integration is with respect to one of the parameters. It is analogous to Theorem 11.5.1 of [15], which is for elements of $D$.

THEOREM 2.10. *Let* $g : \mathbb{R} \to \mathbb{R}$ *be a uniformly continuous function and suppose that* $x_n \to x$ *in* $D_D$ *as* $n \to \infty$ *and* $x \in C_C$. *Then for the integral mapping* $f : D_D \to D_D$ *defined by*

$$f(x)(s,t) = \int_0^t g(x(s,u)) \, du,$$

*we have* $f(x_n) \to f(x)$ *in* $D_D$ *as* $n \to \infty$ *and* $f(x) \in C_C$.

PROOF. First, we show that $f(x) \in C_C$. For each fixed $s \geq 0$, $f(x)(s, \cdot)$ is clearly in $C$. Consider a sequence $(s_n)_{n \geq 1} \subset [0, \infty)$ such that $s_n \to s$ as $n \to \infty$. Then, for each $T > 0$, we have

$$\|f(x)(s_n, \cdot) - f(x)(s, \cdot)\|_T = \left\| \int_0^{\cdot} g(x(s_n, u)) - g(x(s, u)) \, du \right\|_T$$
$$\leq \int_0^T |g(x(s_n, u)) - g(x(s, u))| \, du.$$

The last expression goes to 0 as $n \to \infty$ by bounded convergence and continuity of $g$. Since $f(x) \in C_C$, by Theorem 2.1 it is enough to show that $f(x_n) \to f(x)$ in $(D_D, U, U)$. For each $S, T > 0$, we have

$$(2.7) \qquad \|f(x_n) - f(x)\|_{S,T} \leq \int_0^T \sup_{0 \leq s \leq S} |g(x_n(s,u)) - g(x(s,u))| \, du.$$

Since $x \in C_C$, we have $\|x_n - x\|_{S,T} \to 0$ as $n \to \infty$ for all $S, T > 0$. Thus, by the uniform continuity of $g$, the right-hand side of (2.7) goes to 0 as $n \to \infty$. □

We now give continuity results for regulator mappings that are generalizations of the regulator mapping of Theorem 4.1 of [9].

THEOREM 2.11 (Continuity of regulator map in $D$).   *Consider the integral representation*

$$x(t) = y(t) + \int_0^t h(x(u), u) \, du, \qquad t \geq 0,$$

*where $h : \mathbb{R} \times [0, \infty) \to \mathbb{R}$ satisfies the following properties:*

1. $h(0, t) = 0$ *for all $t \geq 0$.*
2. *For each $T > 0$, there exist $c_1, c_2 > 0$ such that, for each $x_1, x_2 \in D$ and homeomorphism $\lambda$ on $[0, T]$ with strictly positive derivative $\lambda'$, we have*

(2.8)
$$\int_0^t |h(x_1(u), u) - h(x_2(\lambda(u)), \lambda(u))| \, du$$

$$\leq c_1 \|\lambda - e\|_T + c_2 \int_0^t |x_1(u) - x_2(\lambda(u))| \, du$$

*for all $t \in [0, T]$.*

*This integral representation has a unique solution $x$, so that it defines a function $f : D \to D$ mapping $y$ into $x \equiv f(y)$. The function $f$ is continuous provided its domain and range are both equipped with either:* (i) *the topology $U$ or* (ii) *the topology $J_1$. Furthermore, if $y \in C$, then $f(y) \in C$.*

PROOF.   That $f$ is a well-defined map from $D$ to $D$ follows as in the proof of Theorem 4.1 of [9]. To prove continuity in the topology of uniform convergence over bounded intervals, notice that the condition (2.8) implies that there exists $c_2 > 0$ such that for any $x_1, x_2 \in D$ we have

$$\int_0^t |h(x_1(u), u) - h(x_2(u), u)| \, du \leq c_2 \int_0^t |x_1(u) - x_2(u)| \, du$$

for all $t \geq 0$. The proof then follows using Gronwall's inequality as in the proof of Theorem 4.1 of [9].

We now show that $f$ is continuous when both its domain and range are equipped with the $J_1$ topology. We mostly follow the proof of Theorem 4.1 of [9]. Let $T > 0$ be a continuity point of $y \in D$ and suppose $y_n \to y$ in $D([0, T], \mathbb{R})$ with the $J_1$ topology. Then there exists a sequence of increasing homeomorphisms $(\lambda_n)_{n \geq 1}$ of the interval $[0, T]$ such that $\|\lambda_n - e\|_t \to 0$ and $\|y_n - y \circ \lambda_n\|_T \to 0$ as $n \to \infty$. It suffices to consider homeomorphisms $(\lambda_n)_{n \geq 1}$ that are absolutely continuous with respect to Lebesgue measure on $[0, T]$ having derivatives $\lambda_n'$ satisfying $\|\lambda_n' - 1\|_T \to 0$ as $n \to \infty$. Then for $N$ large enough so that $\inf_{0 \leq u \leq T} \lambda_n'(u) > 0$ for all $n > N$ we have

$$|x_n(t) - x(\lambda_n(t))|$$

$$\leq |y_n(t) - y(\lambda_n(t))| + \left| \int_0^t h(x_n(u), u) \, du - \int_0^{\lambda_n(t)} h(x(u), u) \, du \right|$$

$$\leq |y_n(t) - y(\lambda_n(t))|$$

$$+ \left| \int_0^t h(x_n(u), u)\, du - \int_0^y h(x(\lambda_n(u)), \lambda_n(u))\lambda_n'(u)\, du \right|$$

$$\leq |y_n(t) - y(\lambda_n(t))| + \|\lambda_n' - 1\|_T \int_0^T |h(x(\lambda_n(u)), \lambda_n(u))|\, du$$

$$+ \int_0^t |h(x_n(u), u) - h(x(\lambda_n(u)), \lambda_n(u))|\, du$$

$$\leq |y_n(t) - y(\lambda_n(t))| + \|\lambda_n' - 1\|_T c_1 \|x\|_T T + c_1 \|\lambda_n - e\|_T$$

$$+ c_2 \int_0^t |x_n(u) - x(\lambda_n(u))|\, du$$

for all $n > N$ and $0 \leq t \leq T$, where the last inequality follows from (2.8). Since $\|y_n - y \circ \lambda_n\|_T$, $\|\lambda_n' - 1\|_T$ and $\|\lambda_n - e\|_T$ converge to 0 as $n \to \infty$, the result follows by an application of Gronwall's inequality.

Finally, the inheritance of continuity is straightforward as in the proof of Theorem 4.1 of [9]. $\quad\square$

Note that the continuity of $f\colon (D, J_1) \to (D, J_1)$ established in Theorem 2.11 implies measurability using the Kolmogorov $\sigma$-field. We now extend Theorem 2.11 to $D_D$.

THEOREM 2.12 (Continuity of regulator map in $D_D$).    *Consider the integral representation*

$$x(s, t) = y(s, t) + \int_0^t h(x(s, u), s, u)\, du, \qquad t \geq 0,$$

*where $h\colon \mathbb{R} \times [0, \infty)^2 \to \mathbb{R}$ satisfies the following properties:*

1. *$h(0, s, t) = 0$ for all $s, t \geq 0$.*
2. *For each $T > 0$, there exist $c_1, c_2 > 0$ such that for each $s_1, s_2 \geq 0$, $x_1, x_2 \in D_D$, and homeomorphism $\lambda$ on $[0, T]$ with strictly positive derivative $\lambda'$ we have*

$$\int_0^t |h(x_1(s_1, u), s_1, u) - h(x_2(s_2, \lambda(u)), s_2, \lambda(u))|\, du$$

$$\leq c_1 |s_1 - s_2| + c_2 \|\lambda - e\|_T + c_3 \int_0^t |x_1(s_1, u) - x_2(s_2, \lambda(u))|\, du$$

*for all $t \in [0, T]$. This integral representation has a unique solution $x$, so that it defines a function $f\colon D_D \to D_D$ mapping $y$ into $x \equiv f(y)$. The function $f$ is measurable provided its domain and range are equipped with the topology $(D_D, J_1, J_1)$ and continuous provided its domain and range are equipped with the topology $(D_D, U, U)$. Furthermore, if $y \in C_C$, then $f(y) \in C_C$.*

PROOF. Consider our integral representation for each fixed $s \geq 0$. In other words, consider the function $f^s \colon D_D \to D$ mapping $y \in D_D$ to the solution $x \in D$ of

$$x(t) = y(s,t) + \int_0^t h(x(u), s, u) \, du.$$

By Theorem 2.11, $f^s$ is a well-defined function that is continuous in the topology of uniform convergence over bounded intervals. Therefore, our two-parameter integral representation has a unique solution since it has a unique solution for each fixed $s \geq 0$.

To show that for $y \in D_D$, $x \equiv f(y)$ is also in $D_D$, consider a sequence $(s_n)_{n \geq 1} \subseteq [0, \infty)$ such that $s_n \downarrow s$ as $n \to \infty$. Since $y \in D_D$, there exists a family of homeomorphisms $(\lambda_n)_{n \geq 1}$ on $[0, T]$ such that $\|\lambda_n - e\|_T \to 0$ and $\|y^{s_n} - y^s \circ \lambda_n\|_T \to 0$ as $n \to \infty$. Again, it suffices to consider homeomorphisms $(\lambda_n)_{n \geq 1}$ that are absolutely continuous with respect to Lebesgue measure on $[0, T]$ having derivatives $\lambda_n'$ satisfying $\|\lambda_n' - 1\|_T \to 0$ as $n \to \infty$. Then for $N$ large enough so that $\inf_{0 \leq u \leq T} \lambda_n'(u) > 0$ for all $n > N$ we have

$$|x(s_n, t) - x(s, \lambda_n(t))|$$
$$\leq |y(s_n, t) - y(s, \lambda_n(t))|$$
$$\quad + \left| \int_0^t h(x(s_n, u), s_n, u) \, du - \int_0^{\lambda_n(t)} h(x(s, u), s, u) \, du \right|$$
$$\leq |y(s_n, t) - y(s, \lambda_n(t))|$$
$$\quad + \left| \int_0^t h(x(s_n, u), s_n, u) - h(x(s, \lambda_n(u)), s, \lambda_n(u)) \lambda_n'(u) \, du \right|$$
$$\leq |y(s_n, t) - y(s, \lambda_n(t))|$$
$$\quad + \left| \int_0^t h(x(s_n, u), s_n, u) - h(x(s, \lambda_n(u)), s, \lambda_n(u)) \lambda_n'(u) \, du \right|$$
$$\leq |y(s_n, t) - y(s, \lambda_n(t))| + \|\lambda_n' - 1\|_T \int_0^T |h(x(s, \lambda_n(u)), s, \lambda_n(u))| \, du$$
$$\quad + \int_0^t |h(x(s_n, u), s_n, u) - h(x(s, \lambda_n(u)), s, \lambda_n(u))| \, du$$
$$\leq |y(s_n, t) - y(s, \lambda_n(t))| + \|\lambda_n' - 1\|_T c_2 \|x^s\|_T T + c_1 |s_n - s|$$
$$\quad + c_2 \|\lambda_n - e\|_T + c_3 \int_0^t |x(s_n, u) - x(s, \lambda_n(u))| \, du$$

for $n > N$ and $0 \leq t \leq T$. The last inequality follows from conditions 1 and 2. Since $\|y^{s_n} - y^s \circ \lambda_n\|_T$, $\|\lambda_n' - 1\|_T$, $|s_n - s|$ and $\|\lambda_n - e\|_T$ converge to 0 as $n \to \infty$, the result follows by an application of Gronwall's inequality. That left-hand limits exist for $x = f(y)$ follows using a similar argument.

Since the Borel $\sigma$-field on $D_D$ coincides with the Kolmogorov $\sigma$-field generated by the coordinate projections, in order to show measurability, it suffices to show that the map $f^s$ defined above is measurable for each $s \geq 0$. But $f^s$ is simply the composition of the projection $y \mapsto y^s$ and the regulator map $f_1^s$ on $D$ defined in Theorem 2.11 with the function $h_1^s \colon \mathbb{R} \times [0, \infty) \to \mathbb{R}$ given by $h_1^s(x, u) = h(x, s, u)$, for $x \in \mathbb{R}$, $u \geq 0$. Since both the projection $y \mapsto y^s$ and the regulator map $f_1^s$ are measurable, $f^s$ is measurable. Thus, $f$ is measurable.

Continuity in the topology of uniform convergence over bounded intervals follows using the same argument as in Theorem 4.1 of [9], which uses Gronwall's inequality. To show that $y \in C_C$ implies $x = f(y) \in C_C$ is also straightforward using Gronwall's inequality.  □

2.5. *Projection mapping.*  Finally, we will be needing the following projection lemma.

Lemma 2.3.  *Let $x_n \in D_D$, $x \in C_C$, and suppose $x_n \to x$ in $D_D$ as $n \to \infty$. Define $y_n \in D$ by $y_n(t) \equiv x_n(t, t)$ for all $t \geq 0$. Then $y_n \to y$ in $D$ as $n \to \infty$, where $y(t) \equiv x(t, t)$ for all $t \geq 0$.*

Proof.  Apply Theorem 2.1 and the fact that

$$\|y_n - y\|_T = \sup_{0 \leq t \leq T} |x_n(t, t) - x(t, t)|$$
$$\leq \sup_{0 \leq s, t \leq T} |x_n(s, t) - x(s, t)| = \|x_n - x\|_{T, T} \to 0$$

as $n \to \infty$.  □

3. **Critically loaded.**  In this section and the next, we establish heavy-traffic limits for virtual waiting times for general $G/G/n + G$ queues, assuming we have corresponding heavy-traffic limits for the basic arrival, queue-length, departure, and abandonment processes, and the fluid limit of the queue-length process is identically 1. In Section 5, we apply the main result of this section with the basic $M/M/n + M$ QED result, Theorem 5.1, to obtain a proof of the QED virtual waiting time result for the $M/M/n + M$ model, Theorem 5.2.

There are a number of interesting waiting time performance measures that have been studied in the literature. A customer's potential waiting time is the amount of time the customer spends waiting in queue if the customer were infinitely patient (see [3]). Therefore, if $W$ is the potential waiting time of a particular customer and $X$ is his patience (the amount of time he is willing to wait in the queue before being served), then the customer's actual waiting time is given by $W \wedge X$. Another waiting time performance measure of interest is customer waiting time conditional on being served. An exact steady-state analysis of this performance measure in the Markovian case using Laplace transforms appears in [18], but we do not deal with it here. Our main focus is to study virtual waiting time processes, which are continuous-time process analogs of potential waiting times. A virtual waiting time process is a process $(V(t), t \geq 0)$, where $V(t)$ is the potential waiting time of a hypothetical customer arriving to the queue at time $t$.

We assume we are given a sequence of $G/G/n + G$ models. In model $n \geq 1$, let $X_n \equiv (X_n(t), t \geq 0)$, $A_n \equiv (A_n(t), t \geq 0)$, $D_n \equiv (D_n(t), t \geq 0)$, $L_n \equiv (L_n(t), t \geq 0)$ and $V_n \equiv (V_n(t), t \geq 0)$, where

$X_n(t) \equiv$ number of customers in the system at time $t$,

$A_n(t) \equiv$ number of customers that have arrived to the system by time $t$,

$D_n(t) \equiv$ number of customers that have been served by time $t$,

$L_n(t) \equiv$ number of customers that have abandoned by time $t$,

$V_n(t) \equiv$ potential waiting time of a hypothetical customer

arriving at time $t$.

For each $n \geq 1$, define processes

$$\bar{X}_n \equiv \frac{X_n}{n}, \qquad \bar{A}_n \equiv \frac{A_n}{n}, \qquad \bar{D}_n \equiv \frac{D_n}{n}, \qquad \bar{L}_n \equiv \frac{L_n}{n}, \qquad \bar{V}_n \equiv V_n,$$

and

$$\hat{X}_n \equiv c_n(\bar{X}_n - \bar{X}), \qquad \hat{A}_n \equiv c_n(\bar{A}_n - \bar{A}), \qquad \hat{D}_n \equiv c_n(\bar{D}_n - \bar{D}),$$

$$\hat{L}_n \equiv c_n(\bar{L}_n - \bar{L}), \qquad \hat{V}_n \equiv c_n(\bar{V}_n - \bar{V}),$$

where $\bar{X}, \bar{A}, \bar{D}, \bar{L}, \bar{V} \in C$ and $c_n \to \infty$ and $n/c_n \to \infty$ as $n \to \infty$.

We then have the following result.

THEOREM 3.1. *Consider a sequence of $G/G/n + G$ models and suppose*

$$(3.1) \qquad (\bar{A}_n, \bar{D}_n, \bar{L}_n, \bar{X}_n) \Rightarrow (\bar{A}, \bar{D}, \bar{L}, \bar{X}) \qquad in \; D^4$$

*as $n \to \infty$, where $\bar{A}$, $\bar{D}$ and $\bar{L}$ are continuous deterministic functions, $\bar{D}$ has continuous positive derivative $\bar{D}'$, $\bar{L} = c\eta$ for some $c \geq 0$, and $\bar{X} \equiv 1$. Furthermore, suppose*

$$(3.2) \qquad (\hat{A}_n, \hat{D}_n, \hat{L}_n, \hat{X}_n) \Rightarrow (\hat{A}, \hat{D}, \hat{L}, \hat{X}) \qquad in \; D^4$$

*as $n \to \infty$ where $\hat{A}$, $\hat{D}$, $\hat{L}$ and $\hat{X}$ are continuous processes. Then we have*

$$\hat{V}_n \Rightarrow \frac{\hat{X}^+}{\bar{D}'} \qquad in \; D$$

*as $n \to \infty$ jointly with (3.1) and (3.2).*

PROOF. Recalling the discussion about (1.1) and (1.2), first observe that for each $t \geq 0$ and $n \geq 1$, we can bound the virtual waiting time process from above and below by first-passage times:

$$(3.3) \qquad V_n^l(t) \leq V_n(t) \leq V_n^u(t),$$

where

$$(3.4) \qquad V_n^l(t) \equiv \inf\{s \geq 0 | D_n(t+s) + L_n(t+s)$$
$$\geq X_n(t) + D_n(t) + L_n(t) - (n-1)\},$$

$$(3.5) \qquad V_n^u(t) \equiv \inf\{s \geq 0 | D_n(t+s) \geq X_n(t) + D_n(t) - (n-1)\}.$$

For all $n \geq 1$, define the first-passage-time processes $\bar{Z}_n^l \equiv (\bar{Z}_n^l(t), t \geq 0)$ and $\bar{Z}_n^u \equiv (\bar{Z}_n^u(t), t \geq 0)$, for $n \geq 1$, where

$$(3.6) \qquad \begin{aligned} \bar{Z}_n^l(t) &\equiv \inf\{s \geq 0 | \bar{D}_n(s) + \bar{L}_n(s) \geq \bar{X}_n(t) + \bar{D}_n(t) + \bar{L}_n(t) - (1 - 1/n)\} \\ &= \inf\{s \geq 0 | \bar{D}_n(s) + \bar{L}_n(s) \geq \bar{X}_n(0) + \bar{A}_n(t) - (1 - 1/n)\}, \end{aligned}$$

$$(3.7) \qquad \begin{aligned} \bar{Z}_n^u(t) &\equiv \inf\{s \geq 0 | \bar{D}_n(s) \geq \bar{X}_n(t) + \bar{D}_n(t) - (1 - 1/n)\} \\ &= \inf\{s \geq 0 | \bar{D}_n(s) \geq \bar{X}_n(0) + \bar{A}_n(t) - \bar{L}_n(t) - (1 - 1/n)\} \end{aligned}$$

for $t \geq 0$. For the next step, define the processes $\bar{U}_n^l \equiv \bar{Z}_n^l - e$ and $\bar{U}_n^u \equiv \bar{Z}_n^u - e$ for $n \geq 1$. Then we see that the bounds on our virtual waiting time process (3.4) and (3.5) can be written as $V_n^l = (\bar{U}_n^l)^+$ and $V_n^u = (\bar{U}_n^u)^+$, for $n \geq 1$. We would then want to use the corollary of [11], along with the assumptions that $\bar{X} \equiv 1$ and $\bar{L}' \equiv 0$ to get

$$(3.8) \qquad c_n V_n^l = c_n (\bar{Z}_n^l - e)^+ \Rightarrow \frac{\hat{X}^+}{\bar{D}'},$$

$$(3.9) \qquad c_n V_n^u = c_n (\bar{Z}_n^u - e)^+ \Rightarrow \frac{\hat{X}^+}{\bar{D}'}$$

in $D$ as $n \to \infty$. However, notice that the right-hand side of the condition in (3.7) does not satisfy the conditions of the corollary. In particular, $\bar{X}_n(0) + \bar{A}_n - \bar{L}_n - (1 - 1/n)$ is not necessarily a nondecreasing element of $D$ for each $n \geq 1$. Therefore, we cannot deduce (3.9) directly. We can resolve this problem by using linear interpolations that bound the original processes from above.

Notice from (3.7) that

$$\bar{Z}_n^u = \bar{D}_n^{-1} \circ (\bar{X}_n(0) + \bar{A}_n - \bar{L}_n - (1 - 1/n)) = \bar{D}_n^{-1} \circ \bar{Y}_n,$$

where we define $\bar{Y}_n \equiv (\bar{X}_n(0) + \bar{A}_n - \bar{L}_n - (1 - 1/n))$, for $n \geq 1$. Since $\bar{Y}_n$ is not necessarily a nondecreasing element of $D$, $\bar{Z}_n^u$ is not necessarily an element of $D$ (see Example 13.2.1 of [15]). Therefore, we construct linear interpolations $\tilde{D}_n^{-1}$ and $\tilde{Y}_n$ of $\bar{D}_n^{-1}$ and $\bar{Y}_n$, respectively, in such a way that

$$\bar{D}_n^{-1}(t) \leq \tilde{D}_n^{-1}(t) \quad \text{and} \quad \bar{Y}_n(t) \leq \tilde{Y}_n(t)$$

for all $t \geq 0$, and $\tilde{D}_n^{-1}$ is nondecreasing, for each $n \geq 1$. The construction can be carried out as follows: when stepping down, linearly interpolate between the right endpoint of the previous step and the midpoint of the next step. Similarly, when stepping up, linearly interpolate between the midpoint of the previous step and the left endpoint of the next step. We show how to construct the interpolation in Figure 1.

Then, for each $n \geq 1$, let $\tilde{Z}_n^u = \tilde{D}_n^{-1} \circ \tilde{Y}_n$, and notice that, since $\tilde{D}_n^{-1}$ is nondecreasing, we must have $\bar{Z}_n^u(t) \leq \tilde{Z}_n^u(t)$ for $t \geq 0$ so that

$$(3.10) \qquad V_n^u(t) \leq \tilde{V}_n^u(t) \qquad \text{for } t \geq 0,$$

where $\tilde{V}_n^u \equiv (\tilde{Z}_n^u - e)^+$.

By Lemma 7.1, the error caused by these linear interpolations is asymptotically negligible. Also, by our assumed limit (3.2) and Theorem 13.7.2 of [15],

$$\sqrt{n}(\bar{D}_n^{-1} - \bar{D}^{-1}, \bar{Y}_n - (\bar{A} - \bar{L})) \Rightarrow \left(-\left(\frac{\hat{D}}{\bar{D}'}\right) \circ \bar{D}^{-1}, \hat{A} - \hat{L}\right) \qquad \text{in } D^2$$

as $n \to \infty$, jointly with (3.1) and (3.2). Therefore, using the converging together lemma (Theorem 11.4.7 of [15]) gives us

$$\sqrt{n}(\tilde{D}_n^{-1} - \bar{D}^{-1}, \tilde{Y}_n - (\bar{A} - \bar{L})) \Rightarrow \left(-\left(\frac{\hat{D}}{\bar{D}'}\right) \circ \bar{D}^{-1}, \hat{A} - \hat{L}\right) \qquad \text{in } D^2$$

as $n \to \infty$, jointly with (3.1) and (3.2). Since the limit process in the last limit belongs to $C^2$, the limit holds in the stronger uniform topology. But since the prelimit processes also belong to $C^2$ the limit holds in $C^2$. Using a version of the composition result Theorem 13.3.3 of [15] for elements of $C$ with the continuous mapping theorem gives us

$$c_n \tilde{V}_n^u \equiv c_n(\tilde{Z}_n^u - e)^+ = c_n(\tilde{D}_n^{-1} \circ \tilde{Y}_n - e)^+ \Rightarrow \frac{\hat{X}^+}{\bar{D}'} \qquad \text{in } C$$

as $n \to \infty$. Then, using the bounds (3.10) and (3.3), gives us our result. $\quad\square$

A similar bounding argument does not work in the ED case, because the limits of the bounding processes do not coincide.

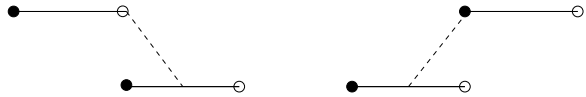

FIG. 1.  *Upper bound linear interpolation for step functions.*

**4. Limits from stopped arrival processes.** In this section, we prove heavy-traffic limits for virtual waiting-times for $G/G/n + G$ queues from associated heavy-traffic limits for basic queueing processes with the arrival process stopped at some time $\tau \geq 0$. The results here will be applied to the Markovian $M/M/n + M$ model in the ED limiting regime in Section 6.1.

We present results for one-dimensional virtual waiting times and virtual-waiting-time processes, respectively, in the next two subsections.

4.1. *One-dimensional virtual waiting times.* To get one-dimensional limits for the virtual waiting time at some fixed time $\tau \geq 0$, the idea is to stop the arrival process to the system at time $\tau$. We add the superscript $\tau$ to all processes defined in Section 3 to denote the same processes but with the arrival process to the system stopped at time $\tau$. We then have for $n \geq 1$,

$$V_n^\tau(t) = \inf\{s \geq 0 | D_n^\tau(t+s) + L_n^\tau(t+s) \geq X_n^\tau(0) + A_n^\tau(t) - (n-1)\}$$

for $\tau, t \geq 0$, and $V_n(t) = V_n^t(t)$ for $t \geq 0$.

For each $\tau \geq 0$, we assume we have the limits

$$(4.1) \qquad (\bar{A}_n^\tau, \bar{D}_n^\tau, \bar{L}_n^\tau, \bar{X}_n^\tau) \Rightarrow (\bar{A}^\tau, \bar{D}^\tau, \bar{L}^\tau, \bar{X}^\tau) \qquad \text{in } D^4$$

as $n \to \infty$, where $\bar{A}^\tau$, $\bar{D}^\tau$, $\bar{L}^\tau$ and $\bar{X}^\tau$ are continuous deterministic functions and $\bar{D}^\tau + \bar{L}^\tau$ has continuous positive derivative $(\bar{D}^\tau + \bar{L}^\tau)'$. In addition, we assume that

$$(4.2) \qquad (\hat{A}_n^\tau, \hat{D}_n^\tau, \hat{L}_n^\tau, \hat{X}_n^\tau) \Rightarrow (\hat{A}^\tau, \hat{D}^\tau, \hat{L}^\tau, \hat{X}^\tau) \qquad \text{in } D^4$$

as $n \to \infty$, where $\hat{A}^\tau$, $\hat{D}^\tau$, $\hat{L}^\tau$ and $\hat{X}^\tau$ are continuous processes.

Then if we define for each $\tau \geq 0$ the first-passage time processes $\bar{Z}_n^\tau \equiv (\bar{Z}_n^\tau(t), t \geq 0)$ for $n \geq 1$, and $\bar{Z}^\tau \equiv (\bar{Z}^\tau(t), t \geq 0)$, where

$$\bar{Z}_n^\tau(t) \equiv \inf\{s \geq 0 | \bar{D}_n^\tau(s) + \bar{L}_n^\tau(s) \geq \bar{X}_n^\tau(0) + \bar{A}_n^\tau(t) - (1 - 1/n)\},$$

$$\bar{Z}^\tau(t) \equiv \inf\{s \geq 0 | \bar{D}^\tau(s) + \bar{L}^\tau(s) \geq \bar{X}^\tau(0) + \bar{A}^\tau(t) - 1\}$$

for $t \geq 0$, the processes $\bar{U}_n^\tau \equiv \bar{Z}_n^\tau - e$ for $n \geq 1$, and $\bar{U}^\tau \equiv \bar{Z}^\tau - e$, then we see that for each $\tau \geq 0$ and $n \geq 1$ our virtual-waiting-time processes can be written as $V_n^\tau = (\bar{U}_n^\tau)^+$. We first prove a limit for the processes $\bar{U}_n^\tau$ and then use this limit to prove our virtual-waiting-time heavy-traffic limits. For each $\tau \geq 0$, define the process $\hat{U}^\tau \equiv (\hat{U}^\tau(t), t \geq 0)$ by

$$(4.3) \qquad \hat{U}^\tau(t) \equiv \frac{\hat{X}_n(0) + \hat{A}^\tau(t) - \hat{D}^\tau(Z^\tau(t)) - \hat{L}^\tau(Z^\tau(t))}{(\bar{D}^\tau + \bar{L}^\tau)'(\bar{Z}^\tau(t))}.$$

Note that $\hat{U}^\tau$ is a continuous process.

LEMMA 4.1. *For each $\tau \geq 0$, under the assumptions (4.1) and (4.2),*

$$c_n(\bar{U}_n^\tau - \bar{U}^\tau) \Rightarrow \hat{U}^\tau \qquad in\ D$$

*as $n \to \infty$, jointly with the limits (4.1) and (4.2), where $\hat{U}$ is given in (4.3).*

PROOF. Applying the corollary of [11], we get

$$c_n(\bar{U}_n^\tau - \bar{U}^\tau) = c_n[(\bar{Z}_n^\tau - e) - (\bar{Z}^\tau - e)] = c_n(\bar{Z}_n^\tau - \bar{Z}^\tau) \Rightarrow \hat{U}$$

in $D$ as $n \to \infty$. $\square$

THEOREM 4.1. *For each $\tau \geq 0$, under the assumptions (4.1) and (4.2), we have:*

1. *If $\bar{X}(\tau) = 1$, then*

$$\hat{V}_n(\tau) = c_n\bar{U}_n^\tau(\tau)^+ \Rightarrow \hat{U}^\tau(\tau)^+ = \frac{\hat{X}(\tau)^+}{(\bar{D} + \bar{L})'(\tau)} \qquad in\ \mathbb{R}$$

   *as $n \to \infty$, jointly with the limits (4.1) and (4.2).*

2. *If $\bar{X}(\tau) > 1$, then*

$$\hat{V}_n(\tau) = c_n(\bar{U}_n^\tau(\tau)^+ - \bar{U}^\tau(\tau)) \Rightarrow \hat{U}^\tau(\tau) \qquad in\ \mathbb{R}$$

   *as $n \to \infty$, jointly with the limits (4.1) and (4.2), where $\hat{U}^\tau(\tau)$ is given in (4.3).*

PROOF. Part 1 follows immediately from Lemma 4.1 and the fact that $\bar{Z}^\tau(\tau) = \tau$, $\bar{X}^\tau(\tau) = \bar{X}(\tau)$, $\bar{D}^\tau(\tau) = \bar{D}(\tau)$ and $\bar{L}^\tau(\tau) = \bar{L}(\tau)$ for all $\tau \geq 0$.

For part 2, notice that the condition $\bar{X}(\tau) > 1$ implies $\bar{U}^\tau(\tau) > 0$. By Lemma 4.1, we get the limit

$$\bar{U}_n^\tau(\tau) \Rightarrow \bar{U}^\tau(\tau) \qquad in\ \mathbb{R}$$

as $n \to \infty$. Therefore, for each $\epsilon > 0$ and $0 < \delta < \bar{U}^\tau(\tau)$, we have

$$\begin{aligned}
&\mathbb{P}[|c_n(\bar{U}_n^\tau(\tau)^+ - \bar{U}^\tau(\tau)) - c_n(\bar{U}_n^\tau(\tau) - \bar{U}^\tau(\tau))| < \epsilon] \\
&\quad = \mathbb{P}[c_n|\bar{U}_n^\tau(\tau)^+ - \bar{U}_n^\tau(\tau)| < \epsilon] \\
&\quad \geq \mathbb{P}[\bar{U}_n^\tau(\tau) > 0] \\
&\quad \geq \mathbb{P}[\bar{U}_n^\tau(\tau) > \bar{U}^\tau(\tau) - \delta] \\
&\quad \geq \mathbb{P}[|\bar{U}_n^\tau(\tau) - \bar{U}^\tau(\tau)| < \delta] \to 1
\end{aligned}$$

(4.4)

as $n \to \infty$. Combining this limit with the convergence in Lemma 4.1 and Theorem 4.1 of [2], we complete the proof. $\square$

4.2. *Virtual waiting time process.* To get an associated limit for the virtual-waiting-time processes, we again assume that the limits (4.1) and (4.2) hold, but where the queueing processes are understood to be processes in $D_D$ indexed by both $\tau$ and $t$. To be precise, for each $n \geq 1$ we define the processes

$$X_n^{(2)} \equiv (X_n^{(2)}(\tau, t), \tau, t \geq 0), \qquad A_n^{(2)} \equiv (A_n^{(2)}(\tau, t), \tau, t \geq 0),$$
$$D_n^{(2)} \equiv (D_n^{(2)}(\tau, t), \tau, t \geq 0), \qquad L_n^{(2)} \equiv (L_n^{(2)}(\tau, t), \tau, t \geq 0),$$

so that

$$X_n^{(2)}(\tau, t) \equiv X_n^{\tau}(t), \qquad A_n^{(2)}(\tau, t) \equiv A_n^{\tau}(t),$$
$$D_n^{(2)}(\tau, t) \equiv D_n^{\tau}(t), \qquad L_n^{(2)}(\tau, t) \equiv L_n^{\tau}(t)$$

for $\tau, t \geq 0$ and $n \geq 1$. Next, define the scaled processes

$$\bar{X}_n^{(2)} \equiv \frac{X_n^{(2)}}{n}, \qquad \bar{A}_n^{(2)} \equiv \frac{A_n^{(2)}}{n}, \qquad \bar{D}_n^{(2)} \equiv \frac{D_n^{(2)}}{n}, \qquad \bar{L}_n^{(2)} \equiv \frac{L_n^{(2)}}{n},$$

and

$$\hat{X}_n^{(2)} \equiv c_n(\bar{X}_n^{(2)} - \bar{X}^{(2)}), \qquad \hat{A}_n^{(2)} \equiv c_n(\bar{A}_n^{(2)} - \bar{A}^{(2)}),$$
$$\hat{D}_n^{(2)} \equiv c_n(\bar{D}_n^{(2)} - \bar{D}^{(2)}), \qquad \hat{L}_n^{(2)} \equiv c_n(\bar{L}_n^{(2)} - \bar{L}^{(2)}),$$

where $c_n \to \infty$ and $n/c_n \to \infty$ as $n \to \infty$, and $\bar{X}^{(2)}$, $\bar{A}^{(2)}$, $\bar{D}^{(2)}$ and $\bar{L}^{(2)}$ are given by

$$(4.5) \qquad \begin{aligned} \bar{X}^{(2)}(\tau, t) &\equiv \bar{X}^{\tau}(t), \qquad \bar{A}^{(2)}(\tau, t) \equiv \bar{A}^{\tau}(t), \\ \bar{D}^{(2)}(\tau, t) &\equiv \bar{D}^{\tau}(t), \qquad \bar{L}^{(2)}(\tau, t) \equiv \bar{L}^{\tau}(t). \end{aligned}$$

Note that $\bar{A}^{(2)}$, $\bar{D}^{(2)}$, $\bar{L}^{(2)}$ and $\bar{X}^{(2)}$ are continuous deterministic functions and, for each $\tau \geq 0$, $\bar{D}^{(2)}(\tau, \cdot) + \bar{L}^{(2)}(\tau, \cdot)$ has continuous positive derivative. Then we assume we have the limits

$$(4.6) \qquad (\bar{A}_n^{(2)}, \bar{D}_n^{(2)}, \bar{L}_n^{(2)}, \bar{X}_n^{(2)}) \Rightarrow (\bar{A}^{(2)}, \bar{D}^{(2)}, \bar{L}^{(2)}, \bar{X}^{(2)}) \qquad \text{in } D_D^4$$

as $n \to \infty$, and

$$(4.7) \qquad (\hat{A}_n^{(2)}, \hat{D}_n^{(2)}, \hat{L}_n^{(2)}, \hat{X}_n^{(2)}) \Rightarrow (\hat{A}^{(2)}, \hat{D}^{(2)}, \hat{L}^{(2)}, \hat{X}^{(2)}) \qquad \text{in } D_D^4$$

as $n \to \infty$, where $\hat{X}^{(2)}$, $\hat{A}^{(2)}$, $\hat{D}^{(2)}$ and $\hat{L}^{(2)}$ are given by

$$(4.8) \qquad \begin{aligned} \hat{X}^{(2)}(\tau, t) &\equiv \hat{X}^{\tau}(t), \qquad \hat{A}^{(2)}(\tau, t) \equiv \hat{A}^{\tau}(t), \\ \hat{D}^{(2)}(\tau, t) &\equiv \hat{D}^{\tau}(t), \qquad \hat{L}^{(2)}(\tau, t) \equiv \hat{L}^{\tau}(t). \end{aligned}$$

Note that $\hat{X}^{(2)}, \hat{A}^{(2)}, \hat{D}^{(2)}, \hat{L}^{(2)} \in C_C$.

Next, define the processes $\bar{U}_n^{(2)} \equiv (\bar{U}_n^{(2)}(\tau,t), \tau, t \geq 0)$, for $n \geq 1$, and the process $\hat{U}^{(2)} \equiv (\hat{U}^{(2)}(\tau,t), \tau, t \geq 0)$ so that

$$(4.9) \qquad \bar{U}_n^{(2)}(\tau,t) \equiv \bar{U}_n^{\tau}(t), \qquad \hat{U}^{(2)}(\tau,t) \equiv \hat{U}^{\tau}(t)$$

for $\tau, t \geq 0$ and $n \geq 1$. Then, as in Lemma 4.1, we have

LEMMA 4.2. *Under the assumptions (4.6) and (4.7), we have*

$$(4.10) \qquad c_n(\bar{U}_n^{(2)} - \bar{U}^{(2)}) \Rightarrow \hat{U}^{(2)} \qquad in\ D_D$$

*as $n \to \infty$ jointly with the limits (4.6) and (4.7), where $\bar{U}^{(2)}$ and $\hat{U}^{(2)}$ are given by (4.9).*

PROOF. The proof follows the proof of Lemma 4.1, except instead of using Corollary 1 of [11] we use Theorem 2.9 here. □

Now define the projections $\bar{U}_n^0 \equiv (\bar{U}_n^0(t), t \geq 0)$, for $n \geq 1$, and $\hat{U}^0 \equiv (\hat{U}^0(t), t \geq 0)$ so that

$$\bar{U}_n^0(t) \equiv \bar{U}_n^{(2)}(t,t), \qquad \hat{U}_n^0(t) \equiv \hat{U}_n^{(2)}(t,t)$$

for $t \geq 0$. Then we have

THEOREM 4.2. *Under the assumptions (4.6), (4.7) and $\hat{X}^{(2)} \in C_C$, we have the following:*

1. *If $\bar{X}(\tau) = 1$ for all $\tau \geq 0$, then*

$$\hat{V}_n = c_n(\bar{U}_n^0)^+ \Rightarrow (\hat{U}^0)^+ = \frac{\hat{X}^+}{(\bar{D} + \bar{L})'} \qquad in\ D$$

   *as $n \to \infty$, jointly with the limits (4.1) and (4.2).*
2. *If $\inf_{\tau \geq 0} \bar{X}(\tau) > 1$, then*

$$\hat{V}_n = c_n((\bar{U}_n^0)^+ - \bar{U}^0) \Rightarrow \hat{U}^0 \qquad in\ D$$

   *as $n \to \infty$, jointly with the limits (4.6) and (4.7).*

PROOF. The limit in Part 1 follows immediately from the limit (4.10) by an application of the continuous mapping theorem with the projection map of Lemma 2.3. The second equality then follows from the fact that $\bar{Z}^t(t) = t$, $\bar{X}^t(t) = \bar{X}(t)$, $\bar{D}^t(t) = \bar{D}(t)$ and $\bar{L}^t(t) = \bar{L}(t)$ for all $t \geq 0$.

For Part 2, notice the condition $\inf_{\tau \geq 0} \bar{X}(\tau) > 1$ implies $\inf_{\tau \geq 0} \bar{U}^\tau(\tau) > 0$ so that as in (4.4) we get for all $\epsilon > 0$ and $T > 0$,

$$\mathbb{P}[\|c_n((\bar{U}_n^0)^+ - \bar{U}^0) - c_n(\bar{U}_n^0 - \bar{U}^0)\|_T < \epsilon] \to 1 \qquad as\ n \to \infty.$$

Combining this limit with the limit (4.10) and the continuous mapping theorem with the projection map of Lemma 2.3 completes the proof. □

**5. QED.** In this section and the next, we apply the results of the previous two sections to prove heavy-traffic limits for virtual waiting times in the Markovian $M/M/n + M$ model in the QED and ED limiting regimes, respectively. We assume that interarrival times, service times and patience times after time 0 are mutually independent and independent of the initial conditions. To construct our heavy-traffic limits, we consider a sequence of $M/M/n + M$ models indexed by $n \geq 1$. The $n$th model has $n$ servers, each with service rate $\mu$, arrival rate $\lambda_n$, and individual abandonment rate $\theta$. The QED limiting regime of [5] is characterized by the limit

$$\lim_{n \to \infty} \sqrt{n}(1 - \rho_n) = \beta, \qquad -\infty < \beta < \infty,$$

where $\rho_n \equiv \lambda_n/n\mu$ is the traffic intensity in the $n$th model. Let the basic processes associated with the models, i.e., $X_n$, $A_n$, $D_n$, $L_n$, $V_n$, $\bar{X}_n$, $\bar{A}_n$, $\bar{D}_n$, $\bar{L}_n$, $\bar{V}_n$, $\hat{X}_n$, $\hat{A}_n$, $\hat{D}_n$, $\hat{L}_n$ and $\hat{V}_n$, for $n \geq 1$, be defined as in Section 3. For this QED case, under initial conditions such that $\bar{X}_n(0) \Rightarrow 1$ in $\mathbb{R}$ as $n \to \infty$, we will have the fluid limits

$$(5.1) \qquad \bar{X} \equiv \eta, \qquad \bar{A} \equiv \mu e, \qquad \bar{D} \equiv \mu e, \qquad \bar{L} \equiv 0, \qquad \bar{V} \equiv 0.$$

Since $\bar{X} \equiv \eta$ and $\bar{L}' \equiv 0$, the results of Section 3 apply.

We start with the following theorem (Theorem 2 of [3]) for all processes except for $V_n$. See Theorem 7.1 of [9] for an alternative proof.

THEOREM 5.1. *Consider the sequence of $M/M/n + M$ models defined above. If $\bar{X}_n(0) \Rightarrow \bar{X}(0)$ in $\mathbb{R}$ as $n \to \infty$, then we have the fluid limit*

$$(5.2) \qquad (\bar{A}_n, \bar{D}_n, \bar{L}_n, \bar{X}_n) \Rightarrow (\bar{A}, \bar{D}, \bar{L}, \bar{X}) \qquad in \ D^4$$

*as $n \to \infty$, where the limit is given in (5.1). Furthermore, if $\hat{X}_n(0) \Rightarrow \hat{X}(0)$ in $\mathbb{R}$ as $n \to \infty$, then we have the diffusion limit*

$$(\hat{A}_n, \hat{D}_n, \hat{L}_n, \hat{X}_n) \Rightarrow (\hat{A}, \hat{D}, \hat{L}, \hat{X}) \qquad in \ D^4$$

*as $n \to \infty$, where $\hat{X}$ is the unique solution to the stochastic differential equation*

$$d\hat{X}(t) = -\mu\beta - \mu(\hat{X}(t) \wedge 0) \, dt - \theta(\hat{X}(t) \vee 0) \, dt + \sqrt{\mu} \, dB_1(t) - \sqrt{\mu} \, dB_2(t),$$

*and*

$$\hat{A}(t) \equiv \sqrt{\mu} B_1(t),$$

$$\hat{D}(t) \equiv \sqrt{\mu} B_2(t) + \mu \int_0^t (\hat{X}(s) \wedge 0) \, ds,$$

$$\hat{L}(t) \equiv \theta \int_0^t (\hat{X}(s) \vee 0) \, ds$$

*for $t \geq 0$, where $(B_1, B_2)$ is standard Brownian motion in two dimensions and is independent of $\hat{X}(0)$.*

Then, applying Theorem 3.1, we obtain:

THEOREM 5.2.   *Under the assumptions of Theorem 5.1,*

$$\hat{V}_n \Rightarrow \frac{\hat{X}^+}{\mu} \qquad in \ D$$

*as $n \to \infty$, jointly with the limits in Theorem 5.1.*

We can use Theorem 5.2 to get a heavy-traffic limit for potential waiting times of customers in the queue (see [13]). Let $\hat{W}_n \equiv (\sqrt{n} W_n^{\lfloor nt \rfloor+1}, t \geq 0)$, where $W_n^i$ denotes the potential waiting time of the $i$th customer to enter system $n$, for $i \geq 1$, $n \geq 1$. Then using the fact that

$$W_n^{\lfloor nt \rfloor+1} = V_n(T_n^{\lfloor nt \rfloor+1}), \qquad t \geq 0,$$

where $T_n^i$ is the arrival time of the $i$th customer in the $n$th system, for $i \geq 1$, $n \geq 1$, along with the continuous mapping theorem with the composition map, we get

COROLLARY 5.1.   *Under the assumptions of Theorem 5.1,*

$$\hat{W}_n \Rightarrow \frac{\hat{X}(\cdot/\mu)^+}{\mu} \qquad in \ D$$

*as $n \to \infty$, jointly with the limits in Theorem 5.1.*

Finally, we prove convergence of steady-state virtual waiting times. The argument is given in the proof of Theorem 3 of [3], but we repeat it here for completeness.

COROLLARY 5.2.   *There exist random variables $\hat{X}(\infty)$ and $\hat{V}_n(\infty)$ such that*

$$(5.3) \qquad\qquad \hat{X}(t) \Rightarrow \hat{X}(\infty) \qquad in \ \mathbb{R} \ as \ t \to \infty,$$

$$(5.4) \qquad\qquad \hat{V}_n(t) \Rightarrow \hat{V}_n(\infty) \qquad in \ \mathbb{R} \ as \ t \to \infty$$

*for $n \geq 1$. Furthermore, under the assumptions of Theorem 5.1,*

$$\hat{V}_n(\infty) \Rightarrow \frac{\hat{X}(\infty)^+}{\mu} \qquad in \ \mathbb{R}$$

*as $n \to \infty$.*

PROOF.   We first recall from [3] that $\hat{X}(\infty)$ is indeed well defined. Also, recall from Theorem 2* of [3] that we can interchange the heavy-traffic and steady-state limits of $\hat{X}_n(t)$. We can use this fact to deduce the interchange for $\hat{V}_n(t)$ as follows. Note that each of the process $X_n$, $n \geq 1$, possesses a stationary distribution since it is simply a birth–death process that can easily be seen to be positive recurrent. For each $n \geq 1$, define $X_n(0)$ to have this distribution so that we have a stationary version of the process $X_n$. Because of the birth–death structure of the processes $X_n$, $n \geq 1$, for each $t \geq 0$ and $n \geq 1$, $V_n(t)$ can be represented in terms of $X_n(t)$ by

$$V_n(t) = \sum_{i=1}^{(X_n(t)-n)^+} X_{i,n},$$

where $X_{i,n}$ represents the time between the $(i-1)$th and $i$th departure (either an abandonment or a service completion) from the queue in model $n \geq 1$. Therefore, for each $n \geq 1$, $\hat{V}_n$ must also be a stationary process so that $\hat{V}_n(t)$ has the same distribution for all $0 \leq t < \infty$. Thus, for each $n \geq 1$, the random variable $\hat{V}_n(\infty)$ defined in (5.4) is well defined with the same distribution as $\hat{V}_n(t)$, for any fixed $t \geq 0$. Projecting the result of Theorem 5.2 onto a fixed $t \geq 0$ and using the fact that we can interchange the heavy-traffic and steady-state limits of $\hat{X}_n(t)$ then gives us our result.    $\square$

**6. ED.** We now study waiting-time asymptotics in the ED limiting regime using the results of Section 4. We again have a sequence of $M/M/n + M$ models with $n$ servers in model $n \geq 1$, but now customers arrive to the system at the rate $\lambda_n \equiv n\lambda$, for some $\lambda > 0$. In all the models, each server has service rate $\mu < \lambda$ and each customer has abandonment rate $\theta$. As in Section 5, we assume that interarrival times, service times and patience times after time 0 are mutually independent and independent of the initial conditions. Again, let the basic processes associated with the models, i.e., $X_n$, $A_n$, $D_n$, $L_n$, $V_n$, $\bar{X}_n$, $\bar{A}_n$, $\bar{D}_n$, $\bar{L}_n$, $\bar{V}_n$, $\hat{X}_n$, $\hat{A}_n$, $\hat{D}_n$, $\hat{L}_n$ and $\hat{V}_n$, for $n \geq 1$, be defined as in Section 3. Under initial conditions such that $\bar{X}_n(0) \Rightarrow 1$ in $\mathbb{R}$ as $n \to \infty$, we will have the fluid limits

$$
\begin{aligned}
&\bar{X} \equiv (q+1)\eta, \qquad q \equiv \frac{\lambda - \mu}{\theta}, \qquad \bar{A} \equiv \lambda e, \qquad \bar{D} \equiv \mu e, \\
&\bar{L} \equiv \theta q e = (\lambda - \mu)e, \qquad \bar{V} \equiv w\eta, \qquad w \equiv \frac{1}{\theta}\ln\left(\frac{\lambda}{\mu}\right)
\end{aligned}
$$

(6.1)

for $t \geq 0$.

We start by stating ED fluid and diffusion limits for $A_n$, $D_n$, $L_n$ and $X_n$. These appear in Theorem 2.1 and Corollary 2.1 of [17].

THEOREM 6.1.    *Consider the sequence of $M/M/n + M$ models defined above. If $\bar{X}_n(0) \Rightarrow \bar{X}(0)$ in $\mathbb{R}$ as $n \to \infty$, then we have the fluid limit*

$$(\bar{A}_n, \bar{D}_n, \bar{L}_n, \bar{X}_n) \Rightarrow (\bar{A}, \bar{D}, \bar{L}, \bar{X}) \qquad in \ D^4$$

(6.2)

as $n \to \infty$, where the limit is given in (6.1). Furthermore, if $\hat{X}_n(0) \Rightarrow \hat{X}(0)$ in $\mathbb{R}$ as $n \to \infty$, then we have the diffusion limit

$$(6.3) \qquad (\hat{A}_n, \hat{D}_n, \hat{L}_n, \hat{X}_n) \Rightarrow (\hat{A}, \hat{D}, \hat{L}, \hat{X}) \qquad in \ D^4$$

as $n \to \infty$, where $\hat{X}$ is the unique solution to the stochastic differential equation

$$d\hat{X}(t) = -\theta \hat{X}(t)\,dt + \sqrt{\lambda}\,dB_1(t) - \sqrt{\mu}\,dB_2(t) - \sqrt{\lambda - \mu}\,dB_3(t),$$

and

$$\hat{A}(t) \equiv \sqrt{\lambda} B_1(t),$$
$$\hat{D}(t) \equiv \sqrt{\mu} B_2(t),$$
$$\hat{L}(t) \equiv \sqrt{\lambda - \mu} B_3(t) + \theta \int_0^t \hat{X}(s)\,ds$$

for $t \geq 0$, where $(B_1, B_2, B_3)$ is standard Brownian motion in three dimensions and is independent of $\hat{X}(0)$.

6.1. *Diffusion limit for virtual waiting time.* We now prove a diffusion limit for the virtual waiting-time process in the ED regime, Theorem 6.3. As discussed in the Introduction, our result is a special case of Theorem 4.1 of [7], but they do not provide a complete proof. We first state the following result, Proposition 6.1, which gives us the limits (4.1) and (4.2) in the case of the $M/M/n + M$ model in the ED limiting regime. It is a special case of Theorems 3.1 and 3.2 of [7], which are proven using strong approximation (see also [6]). The result is also a special case of Theorems 2.1 and 2.2 of [12], which are proved using a different approach. We will then combine Proposition 6.1 with Theorem 4.1 to obtain the one-dimensional heavy-traffic limits for virtual waiting times, Theorem 6.2 below. Afterward, we obtain heavy-traffic limits for virtual waiting time processes in Theorem 6.3.

PROPOSITION 6.1. *Consider the sequence of $M/M/n + M$ models defined above and fix $\tau \geq 0$. If $\bar{X}_n(0) \Rightarrow \bar{X}(0)$ in $\mathbb{R}$ as $n \to \infty$, then we have the fluid limit*

$$(\bar{A}_n^\tau, \bar{D}_n^\tau, \bar{L}_n^\tau, \bar{X}_n^\tau) \Rightarrow (\bar{A}^\tau, \bar{D}^\tau, \bar{L}^\tau, \bar{X}^\tau) \qquad in \ D^4$$

*as $n \to \infty$, jointly with the limit (6.2), where*

$$(6.4) \qquad \bar{X}^\tau(t) \equiv \begin{cases} \dfrac{\lambda e^{-\theta(t-\tau)^+} - \mu}{\theta} + 1, & t < \tau + w, \\ e^{-\mu(t-(\tau+w))}, & t \geq \tau + w, \end{cases}$$

$$(6.5) \qquad \bar{A}^\tau(t) \equiv \lambda(t \wedge \tau),$$

$$(6.6) \qquad \bar{D}^\tau(t) \equiv \mu(t \wedge (\tau + w)) + \frac{1}{\mu}(1 - e^{-\mu(t - (\tau + w))^+}),$$

$$(6.7) \qquad \bar{L}^\tau(t) \equiv (\lambda - \mu)(t \wedge \tau) + \frac{\lambda}{\theta}(1 - e^{-\theta[t - \tau]_0^w}) - \mu[t - \tau]_0^w$$

for $t \geq 0$, where $w \equiv (1/\theta)\ln(\lambda/\mu)$. Furthermore, if $\hat{X}_n(0) \Rightarrow \hat{X}(0)$ in $\mathbb{R}$ as $n \to \infty$, then we have the diffusion limit

$$(\hat{A}_n^\tau, \hat{D}_n^\tau, \hat{L}_n^\tau, \hat{X}_n^\tau) \Rightarrow (\hat{A}^\tau, \hat{D}^\tau, \hat{L}^\tau, \hat{X}^\tau) \qquad in\ D^4$$

as $n \to \infty$, jointly with the limit (6.3), where $\hat{X}^\tau$ is the unique solution to the stochastic integral equation

$$
\begin{aligned}
(6.8) \qquad \hat{X}^\tau(t) = {}& \hat{X}^\tau(0) - \int_0^t [\mu \mathbf{1}_{\{s \geq \tau + w\}} \hat{X}^\tau(s) - \theta \mathbf{1}_{\{s < \tau + w\}} \hat{X}^\tau(s)]\, ds \\
& + \sqrt{\lambda} B_1(t \wedge \tau) - \sqrt{\mu} B_2\left(\int_0^t (\bar{X}^\tau(s) \wedge 1)\, ds\right) \\
& - \sqrt{\theta} B_3\left(\int_0^t (\bar{X}^\tau(s) - 1)^+\, ds\right),
\end{aligned}
$$

and

$$(6.9) \quad \hat{A}^\tau(t) \equiv \sqrt{\lambda} B_1(t \wedge \tau),$$

$$(6.10) \quad \hat{D}^\tau(t) \equiv \mu \int_0^t \mathbf{1}_{\{s \geq \tau + w\}} \hat{X}^\tau(s)\, ds + \sqrt{\mu} B_2\left(\int_0^t (\bar{X}^\tau(s) \wedge 1)\, ds\right),$$

$$(6.11) \quad \hat{L}^\tau(t) \equiv \theta \int_0^t \mathbf{1}_{\{s < \tau + w\}} \hat{X}^\tau(s)\, ds + \sqrt{\theta} B_3\left(\int_0^t (\bar{X}^\tau(s) - 1)^+\, ds\right)$$

for $t \geq 0$, where $(B_1, B_2, B_3)$ is standard Brownian motion in three dimensions and is independent of $\hat{X}(0)$.

We omit the proof of Proposition 6.1 because it is a direct consequence of the two-parameter version in Proposition 6.2 below. The distribution of $\hat{X}^\tau$ in (6.8) is complicated because of the deterministic time transformation in the Brownian motions, but existence and uniqueness follows from Theorem 2.11. Combining Proposition 6.1 with part 2 of Theorem 4.1, we immediately obtain the following one-dimensional result.

THEOREM 6.2. *Under the assumptions of Theorem 6.1, for each $\tau \geq 0$,*

$$\hat{V}_n(\tau) \Rightarrow \frac{\hat{X}^\tau(\tau + w)}{\mu} \qquad in\ \mathbb{R}$$

*as $n \to \infty$, jointly with the limits in Proposition 6.1.*

We now establish a two-parameter version of Proposition 6.1.

PROPOSITION 6.2.  *Consider the sequence of $M/M/n + M$ models defined above. If $\bar{X}_n(0) \Rightarrow \bar{X}(0)$ in $\mathbb{R}$ as $n \to \infty$, then we have the fluid limit*

$$(\bar{A}_n^{(2)}, \bar{D}_n^{(2)}, \bar{L}_n^{(2)}, \bar{X}_n^{(2)}) \Rightarrow (\bar{A}^{(2)}, \bar{D}^{(2)}, \bar{L}^{(2)}, \bar{X}^{(2)}) \qquad in\ D_D^4$$

*as $n \to \infty$, jointly with the limit (6.2), where $\bar{X}^{(2)}$, $\bar{A}^{(2)}$, $\bar{D}^{(2)}$ and $\bar{L}^{(2)}$ are given by (4.5) and (6.4)–(6.7). Furthermore, if $\hat{X}_n(0) \Rightarrow \hat{X}(0)$ in $\mathbb{R}$ as $n \to \infty$, then we have the diffusion limit*

$$(\hat{A}_n^{(2)}, \hat{D}_n^{(2)}, \hat{L}_n^{(2)}, \hat{X}_n^{(2)}) \to (\hat{A}^{(2)}, \hat{D}^{(2)}, \hat{L}^{(2)}, \hat{X}^{(2)}) \qquad in\ D_D^4$$

*as $n \to \infty$, jointly with the limit (6.3), where $\hat{X}^{(2)}$, $\hat{A}^{(2)}$, $\hat{D}^{(2)}$ and $\hat{L}^{(2)}$ are given by (4.8) and (6.8)–(6.10). Furthermore, $(\hat{A}^{(2)}, \hat{D}^{(2)}, \hat{L}^{(2)}, \hat{X}^{(2)}) \in C_C^4$.*

PROOF.  The proof follows the same steps as the proof of Theorem 7.1 of [9] except some extra care must be taken to apply the continuous mapping theorem with our integral mapping. First, after centering, the system equation becomes

$$
\begin{aligned}
X_n^{(2)}(\tau, t) = {} & X_n^{(2)}(\tau, 0) + M_{n,1}^{(2)}(\tau, t) - M_{n,2}^{(2)}(\tau, t) - M_{n,3}^{(2)}(\tau, t) \\
& + \lambda_n(t \wedge \tau) - \mu \int_0^t (X_n^{(2)}(\tau, s) \wedge n)\, ds \\
& - \theta \int_0^t (X_n^{(2)}(\tau, s) - n)^+\, ds
\end{aligned}
\tag{6.12}
$$

for $\tau, t \geq 0$, where the processes $M_{n,i}^{(2)} \equiv (M_{n,i}^{(2)}(\tau, t), \tau, t \geq 0)$ are defined for $i = 1, 2, 3$ by

$$M_{n,1}^{(2)}(\tau, t) \equiv A(\lambda_n(\tau \wedge t)) - \lambda_n(\tau \wedge t),$$

$$M_{n,2}^{(2)}(\tau, t) \equiv S\left(\mu \int_0^t (X_n^{(2)}(\tau, s) \wedge n)\, ds\right) - \mu \int_0^t (X_n^{(2)}(\tau, s) \wedge n)\, ds,$$

$$M_{n,3}^{(2)}(\tau, t) \equiv L\left(\theta \int_0^t (X_n^{(2)}(\tau, s) - n)^+\, ds\right) - \theta \int_0^t (X_n^{(2)}(\tau, s) - n)^+\, ds$$

for $\tau, t \geq 0$, and $A$, $S$ and $L$ are the unit Poisson processes generating arrivals, departures and abandonments in the setting of Theorem 6.1.

*Fluid limit.*  Since for each $n \geq 1$ we have $X_n^{(2)}(\tau, s) \leq X_n(s)$ for all $\tau, s \geq 0$, we have as in Lemma 4.5 of [9]

$$\frac{M_{n,i}^{(2)}}{n} \to 0 \qquad in\ D_D\ \text{w.p. 1}\ \ for\ i = 1, 2, 3,$$
$$\tag{6.13}$$

jointly as $n \to \infty$ using the FSLLN for Poisson processes. We only apply the weaker FWLLN consequence.

Then, using the integral map of Theorem 2.12 with the function $h_1 : \mathbb{R} \times [0, \infty)^2 \to \mathbb{R}$ defined by

$$h_1(x, s, t) = -\mathbf{1}_{\{t \geq s+w\}} \mu x - \mathbf{1}_{\{t < s+w\}} \theta(x-1), \qquad x \in \mathbb{R}, s, t \geq 0,$$

$\bar{X}^\tau$ must satisfy the integral equation

$$\bar{X}^{(2)}(\tau, t) = \bar{X}^{(2)}(\tau, 0) + \lambda(\tau \wedge t) - \mu((\tau + w) \wedge t)$$
$$- \int_0^t \mathbf{1}_{\{s \geq \tau + w\}} \mu \bar{X}^{(2)}(\tau, s) + \mathbf{1}_{\{s < \tau + w\}} \theta(\bar{X}^{(2)}(\tau, s) - 1) \, ds,$$

for $\tau, t \geq 0$. Lemma 7.2 shows that $h_1$ satisfies the conditions of Theorem 2.12. It is easy to see by inspection that this integral equation has solution given by (6.4). We can now use (6.13) to get the limits for $\bar{A}_n^{(2)}$, $\bar{D}_n^{(2)}$ and $\bar{L}_n^{(2)}$ jointly.

*Diffusion limit.* Using our two-parameter integral mapping (Theorem 2.10) and two-parameter composition mapping (Theorem 2.4) with the continuous mapping theorem, we get the limits

$$(6.14) \qquad \hat{M}_{n,i}^{(2)} \equiv \frac{M_{n,i}^{(2)}}{\sqrt{n}} \Rightarrow \hat{M}_i^{(2)} \qquad \text{in } D_D \text{ for } i = 1, 2, 3,$$

jointly as $n \to \infty$, where

$$\hat{M}_1^{(2)}(\tau, t) \equiv B_1((\lambda(\tau \wedge t))),$$
$$\hat{M}_2^{(2)}(\tau, t) \equiv B_2\left(\mu \int_0^t (\bar{X}^{(2)}(\tau, s) \wedge 1) \, ds\right),$$
$$\hat{M}_3^{(2)}(\tau, t) \equiv B_3\left(\theta \int_0^t (\bar{X}^{(2)}(\tau, s) - 1)^+ \, ds\right).$$

Notice that $\bar{X}^{(2)} \in C_C$ so that the applications of the continuous mapping theorem with Theorems 2.10 and 2.4 are justified.

By Lemma 7.4 and system equation (6.12), there exists a sequence of processes $\hat{H}_n \equiv (\hat{H}_n(\tau, t), \tau, t \geq 0)$, $n \geq 1$, such that

$$(6.15) \qquad \hat{H}_n \Rightarrow 0 \qquad \text{in } D_D \text{ as } n \to \infty$$

and

$$X_n^{(2)}(\tau, t) = X_n^{(2)}(\tau, 0) + M_{n,1}^{(2)}(\tau, t) - M_{n,2}^{(2)}(\tau, t) - M_{n,3}^{(2)}(\tau, t)$$
$$+ \lambda_n(\tau \wedge t) - \mu((\tau + w) \wedge t)n + \sqrt{n}\hat{H}_n(\tau, t)$$
$$(6.16)$$

$$- \mu \int_0^t \mathbf{1}_{\{s \geq \tau + w\}} X_n^{(2)}(\tau, s) \, ds$$

$$- \theta \int_0^t \mathbf{1}_{\{s < \tau + w\}} (X_n^{(2)}(\tau, s) - n) \, ds.$$

Centering $X_n^{(2)}(\tau, t)$ by $n\bar{X}^{(2)}(\tau, t)$ in (6.16) and dividing by $\sqrt{n}$, we get

$$\hat{X}_n^{(2)}(\tau, t) = \hat{X}_n^{(2)}(\tau, 0) + \hat{M}_{n,1}^{(2)}(\tau, t) - \hat{M}_{n,2}^{(2)}(\tau, t) - \hat{M}_{n,3}^{(2)}(\tau, t)$$

$$(6.17) \qquad + \hat{H}_n(\tau, t) - \mu \int_0^t \mathbf{1}_{\{s \geq \tau + w\}} \hat{X}_n^{(2)}(\tau, s) \, ds$$

$$- \theta \int_0^t \mathbf{1}_{\{s < \tau + w\}} \hat{X}_n^{(2)}(\tau, s) \, ds$$

for $t \geq 0$.

Define the function $h_0 \colon \mathbb{R} \times [0, \infty)^2 \to \mathbb{R}$ by

$$h_0(x, s, t) = -\mathbf{1}_{\{t \geq s + w\}} \mu x - \mathbf{1}_{\{t < s + w\}} \theta x, \qquad x \in \mathbb{R}, s \geq 0.$$

It follows from Lemma 7.2 with $a = 0$ that this function satisfies the conditions of Theorem 2.12. Applying the continuous mapping theorem with the integral map of Theorem 2.12 and the addition map of Theorem 2.2 to the limits (6.14) and (6.15) gives us our limit for $\hat{X}_n^\tau$. The limits for $\hat{A}_n^{(2)}$ and $\hat{L}_n^{(2)}$ follow directly from the limits for $\hat{M}_{n,1}^{(2)}$ and $\hat{M}_{n,3}^{(2)}$. Then, the limit for $\hat{D}_n^{(2)}$ follows from (6.17) and the continuous mapping theorem with addition.

Theorems 2.4, 2.10 and 2.12 also show that the property of being in $C_C$ is preserved with respect to applications of these two-parameter mappings. $\square$

Now define the process $\hat{X}^w \equiv (\hat{X}^w(t), t \geq 0)$ so that

$$\hat{X}^w(t) \equiv \hat{X}^{(2)}(t, t + w).$$

Then, combining Proposition 6.2 with part 2 of Theorem 4.2 and simplifying $\hat{U}^0$, we obtain the following theorem.

THEOREM 6.3. *Under the assumptions of Theorem 6.1, we have*

$$\hat{V}_n \Rightarrow \frac{\hat{X}^w}{\mu} \qquad in \ D$$

*as $n \to \infty$, jointly with the limits in Theorem 6.1.*

Finally, paralleling Corollary 5.1, we get the following corollary.

COROLLARY 6.1.   *Under the assumptions of Theorem 6.1,*

$$\hat{W}_n \Rightarrow \frac{\hat{X}^w(\cdot/\lambda)}{\mu} \qquad in \ D$$

*as $n \to \infty$, jointly with the limits in Theorem 6.1.*

6.2. *Steady-state virtual waiting time.*   We now prove convergence of steady-state virtual waiting times under diffusion scaling in the ED limiting regime. It is interesting that, although we are only proving a one-dimensional limit, our proof uses convergence in $D$.

Just as in Corollary 5.2, we can define random variables $V_n(\infty)$, for $n \geq 1$, in the ED setting. From (2.26) of [17], we know that

$$(6.18) \qquad V_n(\infty) \Rightarrow w \equiv \frac{1}{\theta} \ln\left(\frac{\lambda}{\mu}\right) \qquad in \ \mathbb{R}$$

as $n \to \infty$. This result is obtained there by analyzing an ODE for the fluid limit of the queue-length process with arrivals turned off. We now give a different argument for (6.18) as well as a refinement. Let $Q_n(\infty)$ denote the steady-state queue-length (excluding customers in service) in system $n \geq 1$. As in the proof of Corollary 5.2, because of the birth–death structure of $Q_n$, $V_n(\infty)$ can be represented directly in terms of $Q_n(\infty)$ by

$$(6.19) \qquad V_n(\infty) = \sum_{i=1}^{Q_n(\infty)} X_{i,n},$$

where $X_{i,n}$ represents the time between the $(i-1)$th and $i$th departure (either an abandonment or a service completion) from the queue in model $n \geq 1$. For each $n \geq 1$, the $X_{i,n}$ are independent of each other and of $Q_n(\infty)$, and each $X_{i,n}$ is exponentially distributed with mean $1/(n\mu + i\theta)$.

First, we prove the following lemma.

LEMMA 6.1.   *Let the triangular array $(X_{i,n})$ be defined as in (6.19). Then*

$$\sqrt{n}\left(\sum_{i=0}^{\lfloor n\cdot \rfloor} X_{i,n} - c\right) \Rightarrow B \circ d \qquad in \ D$$

*as $n \to \infty$, where*

$$c(t) \equiv \frac{1}{\theta} \ln\left(1 + \frac{\theta t}{\mu}\right), \qquad d(t) \equiv \frac{t}{\mu(\mu + \theta t)}$$

*for all $t \geq 0$, and $B$ is standard Brownian motion.*

Proof. It can easily be verified that the asymptotic mean and variance of the sum in (6.19) are $c(t)$ and $d(t)/n$, respectively. Also, the triangular array $(X_{i,n})$ satisfies Lyapunov's condition. Therefore, the functional Lindeberg–Feller central limit theorem applies to the partial sum sequence and we get our result; see Theorem 4.1 of [10]. □

We now give our main result of the section. Note that one could erroneously conclude from Theorem 6.3 that the limiting distribution of $\hat{V}_n(\infty)$ has variance $(1/\mu^2)\,\mathrm{Var}(\hat{X}(\infty)) = \lambda/(\mu^2\theta)$. We see here that this is in fact not the case. What is happening is that the numerator of the limit at a fixed $t \geq 0$, $\hat{X}^{(2)}(t, t+w)$, corresponds to the diffusion-scaled number-in-system observed $w$ units of time *after* stopping the arrival process. This delay in observing the number-in-system is crucial to the analysis and is the reason one can not simply use the steady-state distribution of $\hat{X}$.

Theorem 6.4. *Under the assumptions of Theorem 6.1,*

$$\hat{V}_n(\infty) \Rightarrow N\left(0, \frac{1}{\theta\mu}\right) \qquad in \ \mathbb{R}$$

*as $n \to \infty$, where $\hat{V}_n(\infty) \equiv \sqrt{n}(V_n(\infty) - w)$.*

Proof. Let $Q_n \equiv (Q_n(t), t \geq 0)$, where $Q_n(t)$ is the number of customers in the queue (excluding customers in service) at time $t$ in model $n \geq 1$. Note that for each $n$, the process $Q_n$ possesses a stationary distribution since $Q_n = (X_n - n)^+$ and $X_n$ is simply a positive recurrent birth–death process. Redefine $Q_n(0)$ to have this distribution so that we have a stationary version of the process $Q_n$. By Theorem 2.3 of [17], we have $\bar{Q}_n(0) \Rightarrow q$ in $\mathbb{R}$ as $n \to \infty$, so that applying Theorem 6.1 and the continuous mapping theorem gives us the fluid limit for the queue length process $\bar{Q}_n(t) \Rightarrow q\eta$, in $D$ as $n \to \infty$. Joining this convergence with the convergence in Lemma 6.1 using Theorem 11.4.5 of [15] gives us

$$\left(\sqrt{n}\left(\sum_{i=0}^{\lfloor n\cdot\rfloor} X_{i,n} - c\right), \bar{Q}_n\right) \Rightarrow (B \circ d, q\eta) \qquad in \ D^2$$

as $n \to \infty$. Projecting the second component on to a fixed $t \geq 0$ gives us

$$\left(\sqrt{n}\left(\sum_{i=0}^{\lfloor n\cdot\rfloor} X_{i,n} - c\right), \bar{Q}_n(t)\right) \Rightarrow (B \circ d, q) \qquad in \ D \times \mathbb{R}$$

as $n \to \infty$. Applying Proposition 13.2.1 of [15] gives us

$$(6.20) \qquad \sqrt{n}\left(\sum_{i=0}^{Q_n(t)} X_{i,n} - c(\bar{Q}_n(t))\right) \Rightarrow B(d(q)) \qquad in \ \mathbb{R}$$

as $n \to \infty$. Then,

$$\sqrt{n}\left(\sum_{i=0}^{Q_n(t)} X_{i,n} - c(q)\right) = \sqrt{n}\left(\sum_{i=0}^{Q_n(t)} X_{i,n} - c(\bar{Q}_n(t))\right)$$
$$+ \sqrt{n}(c(\bar{Q}_n(t)) - c(q)),$$

where $c(q) = w$, for $w$ in (6.18). To compute the limit of the second term above, we use the Taylor approximation

$$c(\bar{Q}_n(t)) = c(q) + c'(q)(\bar{Q}_n(t) - q) + c''(q)(\bar{Q}_n(t) - q)^2 + \cdots$$
$$= c(q) + c'(q)(\bar{Q}_n(t) - q) + o_{\mathbb{P}}(n^{-1/2-\epsilon}),$$

where $0 < \epsilon < 1/2$. The second equality above follows from the continuous mapping theorem and the fact that

$$n^{1/2+\epsilon}(\bar{Q}_n(t) - q)^2 = n^{-1/2+\epsilon}(\sqrt{n}(\bar{Q}_n(t) - q))^2 \Rightarrow 0 \qquad \text{in } \mathbb{R}$$

as $n \to \infty$. Observing that $c'(q) = 1/\lambda$ gives us

$$(6.21) \qquad \sqrt{n}(c(\bar{Q}_n(t)) - c(q)) = c'(q)\hat{Q}_n(t) + o_{\mathbb{P}}(n^{-\epsilon}) \Rightarrow \frac{1}{\lambda}\hat{Q}(t) \qquad \text{in } \mathbb{R}$$

as $n \to \infty$. Combining the limits (6.20) and (6.21) and using the continuous mapping theorem with addition we get

$$\sqrt{n}\left(\sum_{i=0}^{Q_n(t)} X_{i,n} - c(q)\right) \Rightarrow B(d(q)) + \frac{1}{\lambda}\hat{Q}(t) \qquad \text{in } \mathbb{R}$$

as $n \to \infty$. Finally, since $\hat{V}_n(\infty)$ is distributed as the stationary distribution of the left-hand side above, and the stationary distribution of the limit is normal with mean 0 and variance

$$d(q) + \frac{1}{\lambda^2}\text{Var}(\hat{Q}) = \frac{q}{\mu\lambda} + \frac{1}{\lambda\theta} = \frac{1}{\theta\mu},$$

we get our result.  □

**7. Additional lemmas.**  We use the following lemma to show that the error caused by our linear interpolations in the proof of Theorem 3.1 is asymptotically negligible.

LEMMA 7.1.  *For the linear interpolations $\tilde{D}_n^{-1}$ and $\tilde{Y}_n$ of $\bar{D}_n^{-1}$ and $\bar{Y}_n$, respectively, we have for all $T > 0$,*

$$\sqrt{n}\|(\tilde{D}_n^{-1}, \tilde{Y}_n) - (\bar{D}_n^{-1}, \bar{Y}_n)\|_T \Rightarrow 0 \qquad \text{in } \mathbb{R}$$

*as $n \to \infty$, jointly with the limits in Theorem 5.1.*

Proof. Notice that for each $T > 0$ our construction gives us

(7.1)    $\|\bar{D}_n^{-1} - \tilde{D}_n^{-1}\|_T = J_T(\bar{D}_n^{-1})$   and   $\|\bar{Y}_n - \tilde{Y}\|_T = J_T(\bar{Y}_n)$,

which both converge to 0 as $n \to \infty$ since $\bar{D}_n^{-1}$ and $\bar{Y}_n$ are continuous by (5.2) and Corollary 13.7.3 of [15]. Therefore, by (7.1), we have

$$\sqrt{n}\|(\tilde{D}_n^{-1}, \tilde{Y}_n) - (\bar{D}_n^{-1}, \bar{Y}_n)\|_T$$
$$= \sqrt{n} \sup_{0 \le t \le T} \|(\tilde{D}_n^{-1}(t) - \bar{D}_n^{-1}(t), \tilde{Y}_n(t) - \bar{Y}_n(t))\|_\infty$$
$$= \sqrt{n} \sup_{0 \le t \le T} \max\{|\tilde{D}_n^{-1}(t) - \bar{D}_n^{-1}(t)|, |\tilde{Y}_n(t) - \bar{Y}_n(t)|\}$$
$$= \sqrt{n} \max\{\|\tilde{D}_n^{-1} - \bar{D}_n^{-1}\|_T, \|\tilde{Y}_n - \bar{Y}_n\|_T\}$$
$$= \sqrt{n} \max\{J_T(\bar{D}_n^{-1}), J_T(\bar{Y}_n)\}$$
$$= \max\{J_T(\sqrt{n}(\bar{D}_n^{-1} - \bar{D}^{-1})), J_T(\sqrt{n}(\bar{Y}_n - (\bar{A} - \bar{L})))\} \Rightarrow 0$$

in $\mathbb{R}$ as $n \to \infty$, since the processes $\hat{D}$ and $\hat{A} - \hat{L}$ are continuous. Here, we assume the maximum norm $\|\cdot\|_\infty$ on $\mathbb{R}^2$ without loss of generality.   □

The following lemma shows that the regulator map used in Proposition 6.1 satisfies the conditions of Theorem 2.12.

Lemma 7.2. *For $a \in \mathbb{R}$, the function $h_a : \mathbb{R} \times [0, \infty)^2 \to \mathbb{R}$ defined by*

$$h_a(x, s, t) = -\mathbf{1}_{\{t \ge s+w\}} \mu x - \mathbf{1}_{\{t < s+w\}} \theta(x - a), \qquad x \in \mathbb{R}, s, t \ge 0,$$

*satisfies the conditions of Theorem 2.12.*

Proof. $h_a$ clearly satisfies condition 1 of Theorem 2.12. We now show $h_a$ satisfies condition 2 of Theorem 2.12. Fix $T > 0$, $s_1, s_2 \ge 0$, $x_1, x_2 \in D_D$, and an increasing homeomorphism $\lambda$ on $[0, T]$ having strictly positive derivative $\lambda'$. We then have

$$\int_0^t |h_a(x_1(s_1, u), s_1, u) - h_a(x_2(s_2, \lambda(u)), s_2, \lambda(u))|\, du$$

$$= \int_0^t |-\mathbf{1}_{\{u \ge s_1+w\}} \mu x_1(s_1, u) - \mathbf{1}_{\{u < s_2+w\}} \theta(x_1(s_1, u) - a)$$
$$+ \mathbf{1}_{\{\lambda(u) \ge s_2+w\}} \mu x_2(s_2, \lambda(u)) + \mathbf{1}_{\{\lambda(u) < s_2+w\}} \theta(x_2(s_2, \lambda(u)) - a)|\, du$$

$$\le \mu \int_0^t |\mathbf{1}_{\{u \ge s_1+w\}} x_1(s_1, u) - \mathbf{1}_{\{\lambda(u) \ge s_2+w\}} x_2(s_2, \lambda(u))|\, du$$

$$+ \theta \int_0^t |\mathbf{1}_{\{u < s_1+w\}} x_1(s_1, u) - \mathbf{1}_{\{\lambda(u) < s_2+w\}} x_2(s_2, \lambda(u))|\, du$$

$$+ \theta a \int_0^t |\mathbf{1}_{\{u < s_1 + w\}} - \mathbb{1}_{\{\lambda(u) < s_2 + w\}}| \, du$$

$$\leq \mu \int_0^t |\mathbf{1}_{\{u \geq s_1 + w\}} - \mathbb{1}_{\{\lambda(u) \geq s_2 + w\}}| |x_2(s_2, \lambda(u))| \, du$$

$$+ \mu \int_0^t \mathbf{1}_{\{u \geq s_1 + w\}} |x_1(s_1, u) - x_2(s_2, \lambda(u))| \, du$$

$$+ \theta \int_0^t |\mathbf{1}_{\{u < s_1 + w\}} - \mathbb{1}_{\{\lambda(u) < s_2 + w\}}| |x_2(s_2, \lambda(u))| \, du$$

$$+ \theta \int_0^t \mathbf{1}_{\{u < s_1 + w\}} |x_1(s_1, u) - x_2(s_2, \lambda(u))| \, du$$

$$+ \theta a \int_0^t |\mathbf{1}_{\{u < s_1 + w\}} - \mathbb{1}_{\{\lambda(u) < s_2 + w\}}| \, du$$

$$\leq \mu \|x_2^{s_2}\|_T + \theta(\|x_2^{s_2}\|_T + a)|\lambda^{-1}((s_2 + w) \wedge T) - ((s_1 + w) \wedge T)|$$

$$+ (\mu + \theta) \int_0^t |x_1(s_1, u) - x_2(s_2, \lambda(u))| \, du$$

$$\leq \mu \|x_2^{s_2}\|_T + \theta(\|x_2^{s_2}\|_T + a)$$

$$\times (|\lambda^{-1}((s_2 + w) \wedge T) - \lambda^{-1}((s_1 + w) \wedge T)| + \|\lambda^{-1} - e\|_T)$$

$$+ (\mu + \theta) \int_0^t |x_1(s_1, u) - x_2(s_2, \lambda(u))| \, du$$

$$\leq \mu \|x_2^{s_2}\|_T + \theta(\|x_2^{s_2}\|_T + a)\|(\lambda^{-1})'\|_T (|s_2 - s_1| + \|\lambda - e\|_T)$$

$$+ (\mu + \theta) \int_0^t |x_1(s_1, u) - x_2(s_2, \lambda(u))| \, du,$$

where the last inequality follows from the mean value theorem. Since

$$\|(\lambda^{-1})'\|_T \leq \frac{1}{\inf_{0 \leq s \leq T} \lambda'(s)} < \infty,$$

this completes the verification of condition 2 of Theorem 2.12.    □

We now prove the approximations used in the proof of Proposition 6.2. First, we have the following lemma, which is essentially a two-parameter version of Theorem 10.1 of [6].

LEMMA 7.3.   *In the context of Section 6, if* $\hat{X}_n(0) \Rightarrow \hat{X}(0)$ *in* $\mathbb{R}$ *as* $n \to \infty$, *then for all* $S, T > 0$,

$$\|\bar{X}_n^{(2)} - \bar{X}^{(2)}\|_{S,T} = o_{\mathbb{P}}\left(\frac{1}{n^{1/2 - \epsilon}}\right).$$

PROOF. For each $n \geq 1$, define the process $\bar{Y}_n \equiv (\bar{Y}_n(\tau, t), \tau, t \geq 0)$ by

$$\bar{Y}_n(\tau, t) \equiv |\bar{X}_n(0) - \bar{X}(0)| + \left| \frac{\lambda_n}{n}(t \wedge \tau) - \lambda(t \wedge \tau) \right|$$

$$+ \frac{1}{n}|A(\lambda_n(t \wedge \tau)) - \lambda_n(t \wedge \tau)|$$

$$+ \frac{1}{n}\left| S\left( n\mu \int_0^t (\bar{X}_n^{(2)}(\tau, s) \wedge 1)\, ds \right) - n\mu \int_0^t (\bar{X}_n^{(2)}(\tau, s) \wedge 1)\, ds \right|$$

$$+ \frac{1}{n}\left| R\left( n\mu \int_0^t (\bar{X}_n^{(2)}(\tau, s) - 1)^+ \, ds \right) \right.$$

$$\left. - n\mu \int_0^t (\bar{X}_n^{(2)}(\tau, s) - 1)^+ \, ds \right|.$$

Then we have for all $\epsilon > 0$

$$(7.2) \qquad \|\bar{Y}_n\|_{S,T} = o_{\mathbb{P}}\left( \frac{1}{n^{1/2-\epsilon}} \right),$$

by the assumed limit $\hat{X}_n(0) \Rightarrow \hat{X}(0)$ in $\mathbb{R}$ as $n \to \infty$, the assumption $\lambda_n = n\lambda$, the limits (6.14), and Theorem 2.2. Now by the integral representation (6.12), for each $n \geq 1$ we have

$$|\bar{X}_n^{(2)}(\tau, t) - \bar{X}^{(2)}(\tau, t)|$$

$$\leq \bar{Y}_n(\tau, t) + \left| \int_0^t \mu(\bar{X}_n^{(2)}(\tau, s) \wedge 1) + \theta(\bar{X}_n^{(2)}(\tau, s) - 1)^+ \, ds \right.$$

$$\left. - \int_0^t \mu(\bar{X}^{(2)}(\tau, s) \wedge 1) + \theta(\bar{X}^{(2)}(\tau, s) - 1)^+ \, ds \right|$$

$$\leq \bar{Y}_n(\tau, t) + (\mu + \theta) \int_0^T |\bar{X}_n^{(2)}(\tau, s) - \bar{X}^{(2)}(\tau, s)|\, ds$$

for all $\tau, t \geq 0$, where the second inequality follows from the fact that the function

$$h(s) = -\mu(s \wedge 1) - \theta(s - 1)^+$$

is Lipschitz. This implies

$$\|\bar{X}_n^{(2)} - \bar{X}^{(2)}\|_{S,T} \leq \|\bar{Y}_n\|_{S,T} + (\mu + \theta) \int_0^T \|\bar{X}_n^{(2)} - \bar{X}^{(2)}\|_{S,s}\, ds.$$

Applying Gronwall's inequality then gives us

$$\|\bar{X}_n^{(2)} - \bar{X}^{(2)}\|_{S,T} \leq e^{(\mu + \theta)T} \|\bar{Y}_n\|_{S,T}.$$

The result now follows from (7.2). □

We use the following approximation in (6.16).

LEMMA 7.4. *Under the assumptions of Proposition 6.1, for all $S, T \geq 0$,*

$$\sup_{\substack{0 \leq \tau \leq S \\ 0 \leq t \leq T}} \left| \mu \int_0^t (X_n^{(2)}(\tau, s) \wedge n) \, ds + \theta \int_0^t (X_n^{(2)}(\tau, s) - n)^+ \, ds \right.$$

$$(7.3) \qquad - \mu((\tau + w) \wedge t)n - \int_0^t \mathbf{1}_{\{s \geq \tau + w\}} \mu X_n^{(2)}(\tau, s)$$

$$\left. - \mathbf{1}_{\{s < \tau + w\}} \theta(X_n^{(2)}(\tau, s) - n) \, ds \right| = o_{\mathbb{P}}(\sqrt{n}).$$

PROOF. The key is to show that for each $\tau \geq 0$ and $\epsilon > 0$ the first and last times $X_n^{(2)}(\tau, \cdot)$ hits $n$ are within $o_{\mathbb{P}}(1/n^{1/2-\epsilon})$ of $\tau + w$ (see Figure 2). Fix $T > \tau + w$. By Lemma 7.3, for all $\delta_1, \delta_2 > 0$, there exists $N > 0$ such that for all $n \geq N$

$$(7.4) \qquad \mathbb{P}[A_n^{\delta_1}] > 1 - \delta_2,$$

where $A_n^{\delta_1}$ is the event

$$A_n^{\delta_1} \equiv \{\|\bar{X}_n^{(2)} - \bar{X}^{(2)}\|_{S,T} < \delta_1 n^{-1/2+\epsilon}\}.$$

On $A_n^{\delta_1}$ we have for all $0 \leq \tau \leq S$,

$$(7.5) \qquad \inf\{s \geq 0 | X_n^{(2)}(\tau, s) < n\} \geq B_l^\tau,$$

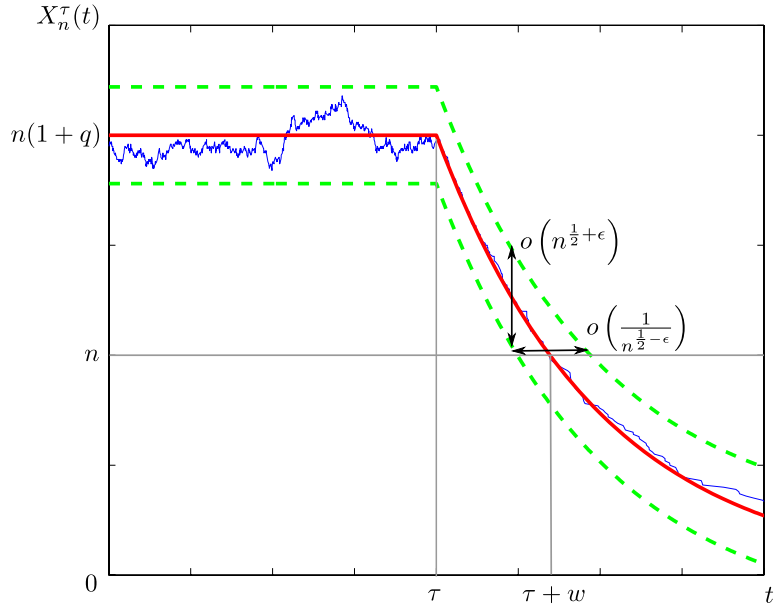

FIG. 2. *The process $X_n^\tau$ deviating from the process $n\bar{X}^\tau$.*

$$(7.6) \qquad \sup\{s \geq 0 | X_n^{(2)}(\tau, s) > n\} \leq B_u^{\tau},$$

where

$$(7.7) \qquad B_l^{\tau} \equiv \inf\left\{s \geq 0 \Big| n\left(\frac{\lambda e^{-\theta(s-\tau)^+} - \mu}{\theta} + 1\right) - \delta_1 n^{1/2+\epsilon} < n\right\},$$

$$(7.8) \qquad B_u^{\tau} \equiv \inf\{s \geq 0 | n e^{-\mu(s-(\tau+w))^+} + \delta_1 n^{1/2+\epsilon} < n\}.$$

It can then be easily checked by manipulating the inequalities in (7.7) and (7.8) that

$$\sup_{0 \leq \tau \leq S} |B_l^{\tau} - (\tau + w)| = o_{\mathbb{P}}\left(\frac{1}{n^{1/2-\epsilon}}\right),$$

$$\sup_{0 \leq \tau \leq S} |B_u^{\tau} - (\tau + w)| = o_{\mathbb{P}}\left(\frac{1}{n^{1/2-\epsilon}}\right).$$

We now show that

$$\sup_{\substack{0 \leq \tau \leq S \\ 0 \leq t \leq T}} \left| \mu \int_0^t (X_n^{(2)}(\tau, s) \wedge n)\, ds \right.$$

$$\left. - \mu((\tau + w) \wedge t)n - \mu \int_0^t \mathbf{1}_{\{s \geq \tau+w\}} X_n^{(2)}(\tau, s)\, ds \right| = o_{\mathbb{P}}(\sqrt{n}).$$

The approximation for the second integral in (7.3) can be shown using the same procedure and this will complete the proof. Let $0 < \epsilon < 1/2$. By (7.5) and (7.6), for all $c > 0$ there exists $N'$ large enough so that for all $n \geq N'$, (7.4) holds and on the event $A_n^{\delta_1}$ we have for all $\tau, t \geq 0$

$$\left| \int_0^t (X_n^{(2)}(\tau, s) \wedge n)\, ds - \left(((\tau + w) \wedge t)n + \int_0^t \mathbf{1}_{\{s \geq \tau+w\}} X_n^{(2)}(\tau, s)\, ds\right) \right|$$

$$\leq \left| \int_0^t (X_n^{(2)}(\tau, s) \wedge n)\, ds \right.$$

$$\left. - \left(\int_0^t \mathbf{1}_{\{s < \tau+w\}} n\, ds + \int_0^t \mathbf{1}_{\{s \geq \tau+w\}} X_n^{(2)}(\tau, s)\, ds\right) \right|$$

$$\leq \int_{\tau+w-c/n^{1/2-\epsilon}}^{\tau+w} (n - X_n^{(2)}(\tau, s))^+\, ds$$

$$+ \int_{\tau+w}^{\tau+w+c/n^{1/2-\epsilon}} (X_n^{(2)}(\tau, s) - n)^+\, ds$$

$$\leq \int_{\tau+w-c/n^{1/2-\epsilon}}^{\tau+w} |X_n^{(2)}(\tau, s) - n\bar{X}^{\tau}(s)|\, ds$$

$$+ \int_{\tau+w}^{\tau+w+c/n^{1/2-\epsilon}} |X_n^{(2)}(\tau, s) - n\bar{X}^{\tau}(s)|\, ds$$

$$\leq c\delta_1 \frac{n^{1/2+\epsilon}}{n^{1/2-\epsilon}} < c\delta_1 n^{1/2}.$$

By (7.4), the above holds with probability greater than $1 - \delta_2$, giving us our result.   $\square$

**Acknowledgments.**   The authors would like to thank Avi Mandelbaum, Bill Massey, Marty Reiman and Sasha Stolyar for sharing their unpublished work and for their support.

Department of Industrial Engineering
and Operations Research
Columbia University
New York, New York 10027-6699
USA
E-mail: rt2146@columbia.edu
        ww2040@columbia.edu